\documentclass [12pt,reqno]{amsart}
\usepackage{times}
\usepackage{amssymb,latexsym,amsmath,amsfonts, eucal}
\usepackage[usenames,dvipsnames,svgnames,table]{xcolor}

\newtheorem{thm}{Theorem}[section]
      \newtheorem{lemma}[thm]{Lemma}
      \newtheorem{prop}[thm]{Proposition}
      
      \newtheorem{defn}[thm]{Definition}
      \newtheorem{examp}[thm]{Example}
      
      \newtheorem{rmk}[thm]{Remark}
      
      \numberwithin{equation}{section}
\title [de Branges spaces based on operator valued reproducing kernels]{ Vector valued de Branges spaces of entire functions based on pairs of  Fredholm operator valued functions and functional model}

\author[ Mahapatra]{Subhankar Mahapatra}

\address{
	Department of Mathematics\\
	Indian Institute of Technology Ropar\\
	140001\\
	India}

\email{ {subhankar.19maz0001@iitrpr.ac.in, subhankarmahapatra95@gmail.com}}

\author[Sarkar]{Santanu Sarkar}

\address{
	Department of Mathematics\\
	Indian Institute of Technology Ropar\\
	140001\\
	India}

\email{ {santanu@iitrpr.ac.in, santanu87@gmail.com.}}

\begin{document}


\keywords{ de Branges spaces with operator valued kernels, multiplication operator, de Branges operator, symmetric operators with infinite deficiency indices, Fredholm operator valued holomorphic functions, Entire operators, Sampling formulas.\\
Accepted: December 02, 2023 Journal of Mathematical Analysis and Applications; doi: https://doi.org/10.1016/j.jmaa.2023.128010.}

\begin{abstract}
\noindent  In this paper, we have considered  vector valued reproducing kernel Hilbert spaces (RKHS) $\mathcal{H}$ of  entire functions associated with operator valued kernel functions. de Branges operators $\mathfrak{E}=(E_-  ,  E_+)$ analogous to de Branges matrices have been constructed with  the help of pairs of Fredholm operator valued entire functions on $\mathfrak{X}$, where $\mathfrak{X}$ is a complex seperable Hilbert space. A few explicit examples of these de Branges operators are also discussed. The newly defined RKHS $\mathcal{B}(\mathfrak{E})$ based on the de Branges operator $\mathfrak{E}=(E_-,E_+)$ has been characterized under some special restrictions. The complete parametrizations and canonical descriptions of all selfadjoint extensions of the closed, symmetric multiplication operator by the independent variable have been given in terms of unitary operators between ranges of reproducing kernels. A sampling formula for the de Branges spaces $\mathcal{B}(\mathfrak{E})$ has been discussed. A particular class of entire operators with infinite deficiency indices has been dealt with and  shown that they can be considered as the multiplication operator for a specific class of these de Branges spaces. Finally, a brief discussion on the connection between the characteristic function of a completely nonunitary contraction operator and the de Branges spaces $\mathcal{B}(\mathfrak{E})$ has been given.
\end{abstract}

\maketitle
\setcounter{tocdepth}{2}

\tableofcontents

\section{Introduction}
\noindent This article aims to develop a general framework of the de Branges theory of vector valued entire functions and establish its connections with the M. G. Krein's theory of a class of entire operators with infinite deficiency indices. The basic theory of reproducing kernel Hilbert spaces (RKHS) was developed by the contribution of many authors (see: \cite{Aronszajn 1950}, \cite{Bergman}, \cite{Mercer 1909}, \cite{Moore 1939}).  It has been observed that the theory of vector valued reproducing kernel Hilbert spaces associated with operator valued kernel functions arise naturally in many areas like probability and stochastic process, machine learning, statistics, etc., and is an active area of research. For example, the articles ( \cite{Micchelli1}, \cite{Micchelli2}) connecting the machine learning theory with the vector valued RKHS theory are motivating. Here we mainly work with RKHS, whose reproducing kernels (RK) are operator valued entire functions.
Throughout this article, $\mathfrak{X}$ is a complex seperable Hilbert space, and $\mathcal{H}$ is an RKHS of $\mathfrak{X}$-valued entire functions. For any $\beta\in\mathbb{C}$, we denote 
 $$\mathcal{H}_\beta=\{g\in\mathcal{H}:g(\beta)=0\}.$$
 Also, by the multiplication operator on $\mathcal{H}$, we mean the operator of multiplication by the independent variable.\\
 In $1959$ with the help of three axioms, L. de Branges introduced the Hilbert spaces of entire functions ($\mathbb{C}$-valued), which were actually RKHS (see \cite{Branges e1}) and are now known as the de Branges spaces. In later years, de Branges published several articles (see, e.g., \cite{Branges e2}, \cite{Branges e3}) and developed the theory as a generalization of Fourier analysis. He also extended his theory for the vector valued analytic functions (see \cite{Branges 1}, \cite{Branges 2}, \cite{Branges 3}) and vector valued entire functions (\cite{Rovnyak}). Over  the time this theory has thrived and made connections with several areas of mathematical analysis, such as the spectral theory of canonical systems, interpolation and sampling. His initial studies can be found in a combined form in the book (\cite{Branges 4}). The  theory of de Branges spaces $\mathcal{B}(\mathfrak{E})$ consisting of $\mathbb{C}^p$-valued entire functions based on a $p\times 2p$ entire matrix valued function $\mathfrak{E}=[E_-~~E_+]$, called the  de Branges matrix, is  also a well developed theory and  appeared to be very crucial for answering direct and inverse problems for canonical systems of differential equations and Dirac-Krein systems . The components of $\mathfrak{E}$ satisfy the following conditions:
\begin{equation}
\det E_+(z)\not\equiv 0 ~\mbox{and}~F:=E_+^{-1}E_-~\mbox{is a}~p\times p~\mbox{inner matrix valued\label{C2} function}.
\end{equation}
The reproducing kernel of $\mathcal{B}(\mathfrak{E})$ is given by
\begin{equation}
K_\xi^{\mathfrak{E}}(z):= \left\{
    \begin{array}{ll}
         \frac{E_+(z)E_+(\xi)^*-E_-(z)E_-(\xi)^*}{\rho_\xi(z)}  & \mbox{if } z \neq \overline{\xi} \\
         \frac{E_+^{'} (\overline{\xi})E_+(\xi)^*- E_-^{'}(\overline{\xi})E_-(\xi)^*}{-2\pi i} & \mbox{if } z = \overline{\xi},
    \end{array} \right. \label{De Branges operator valued kernel}
\end{equation}
which is a $p\times p$ matrix valued entire function. The notation $\rho_\xi(z)$ is clarified later. A comprehensive study of these de Branges spaces with matrix valued reproducing kernels can be found in \cite{ArD08}, \cite{ArD18}. Also  a characterization of the  space $\mathcal{B}(\mathfrak{E})$, analogous  to  problem $50$ of the book \cite{Branges 4}, and it's connection with entire operators having deficiency indices $(p,p)$ are present in \cite{JFA}.\\
 M. G. Krein also studied the Hilbert spaces of entire functions, though his approach was  different from de Branges. He introduced the notion of entire operators and observed the multiplication operator as an operator model of entire operators with arbitrary finite and equal deficiency indices as well as with infinite deficiency indices  in a Hilbert space of entire functions (see: \cite{Krein}, \cite{Krein1}). In recent works for entire operators with $(1,1)$ deficiency indices, these spaces were identified as the de Branges spaces of $\mathbb{C}$-valued entire functions (see \cite{Silva}), and for entire operators with $(p,p)$ deficiency indices ($p$ is arbitrary and finite); these spaces were identified as the de Branges spaces of $\mathbb{C}^p$-valued entire functions (see \cite{JFA}).\\
Our primary goal is to make sense of de Branges operators, likewise the de Branges matrices, i.e., to find a pair of $B(\mathfrak{X})$-valued entire functions $\mathfrak{E}=(E_-,E_+)$ such that the components of $\mathfrak{E}$ would satisfy the conditions like (\ref{C2}) and the kernel function given by (\ref{De Branges operator valued kernel}) should be positive. The undertone of this paper is to observe the transition of the theory of de Branges spaces based on matrix valued reproducing kernels to $B(\mathfrak{X})$-valued reproducing kernels. The two prime operators we shall be considering on $\mathcal{H}$ are the multiplication operator $\mathfrak{T}$ by the independent variable with domain $\mathcal{D}=\{g\in\mathcal{H}:\mathfrak{T}g\in\mathcal{H}\}$ and the generalized backward shift operator defined by
\begin{equation}
    (R_z g)(\xi) := \left\{
    \begin{array}{ll}
         \frac{g(\xi)-g(z)}{\xi-z}  & \mbox{if } \xi \neq z \\
         g'(z) & \mbox{if } \xi = z
    \end{array} \right.
\end{equation}
for every $\xi,z\in\mathbb{C}$. Observing the intimate connections between $\mathcal{H}$, $\mathcal{H}_\beta$, $\mathfrak{T}$ and $R_\beta$ is also one of the main aims of this article.
 \subsection{Plan of the paper}
The paper is organized as follows: The first three sections are dedicated to reviewing the preliminaries for the proposed construction of de Branges operators. Though most of the results in these sections are available in the literature in the matrix setting, to maintain the flow of the study, we mention all the essential results in the operator setting. The main construction of de Branges operators $\mathfrak{E}=(E_-,E_+)$ has been accomplished in section \ref{Se 5}. One of the crucial observations is to consider the components $E_\pm$ of $\mathfrak{E}$ as the Fredholm operator valued entire functions. A particular form of analytic Fredholm theorem (see \cite{Gohberg1}, \cite{Holden}), which we have mentioned in Theorem \ref{Fredholm thm}, turned out to be very crucial for our study. For the basic theory of operator valued holomorphic functions, see \cite{Hille}, \cite{Arendt}, \cite{Gohberg}. In section \ref{Se 6}, we present several examples of de Branges operators. An example of the Fredholm operator valued holomorphic function from the book \cite{Gohberg1} appeared to be very motivating for constructing our example. Section~\ref{Se 7} reviews some results, connecting $\mathcal{H}_\beta$, $\mathfrak{T}$ and $R_\beta$. Also, we discuss the condition for $\mathcal{H}_\beta$ and $\mathcal{H}_{\overline{\beta}}$ to be isometrically isomorphic. Section~\ref{Se 8} discusses a characterization of the newly constructed de Branges spaces $\mathcal{B}(\mathfrak{E})$ corresponding to the de Branges operator $\mathfrak{E}=(E_-,E_+)$, where $E_+(\beta)$ and $E_-(\overline{\beta})$ both are self-adjoint for some $\beta\in\mathbb{C}_+$. This characterization is the vector generalization of problem 50 in \cite{Branges 4} and Theorem $7.1$ in \cite{JFA}.\\
In Section~\ref{Se 9}, under general consideration, we mention the complete parametrizations and canonical descriptions of all selfadjoint extensions of the symmetric multiplication operator $\mathfrak{T}$. Also, a sampling formula for the de Branges spaces has been observed in terms of the eigenfunctions of a selfadjoint extension of $\mathfrak{T}$. The connection between entire operators with infinite deficiency indices and the de Branges spaces $\mathcal{B}(\mathfrak{E})$ is presented in section \ref{Se 10}. Mainly, a particular class of entire operators with infinite deficiency indices have been dealt with and  shown that they can be considered as the multiplication operator for a specific class of de Branges spaces with operator valued RKs. Finally,  in the last section, a brief discussion on the connection between the characteristic function of a completely nonunitary contraction operator and de Branges spaces $\mathcal{B}(\mathfrak{E})$ has been given.

\subsection{Notations} The following notations will be used throughout the paper:\\
 $\mathbb{C}$ the complex plane;  $\mathbb{C}_+$ (resp., $\mathbb{C}_-$) the open upper (resp., lower) half-plane. $\mathbb{D}$ is the open unit disc and $\mathbb{T}$ is the unit circle in the complex plane.\\
 $B(\mathfrak{X})$ denotes the collection of all bounded linear operators on $\mathfrak{X}$.\\
 $\rho_\xi(z)=-2\pi i(z-\overline{\xi})$.\\
 For an operator $A$; $A^*$ denotes the adjoint operator, the notation $A\succ 0$ and $A\succeq 0$ means that $A$ is positive definite and positive semi-definite respectively and $\sigma(A)$ denote the spectrum of $A$.\\ 
 A point $\beta\in\mathbb{C}$ is a point of regular type of $A$ if there exists $d_\beta>0$ such that
 $$||(A-\beta I)g||\geq d_\beta ||g||$$
 for every $g$ in the domain of $A$. $\pi(A)$ denotes the collection of all points of regular type for $A$.\\
 $\ker A$ denotes the kernel of $A$; $\mbox{rng}~ A$ denotes the range of $A$, and $\overline{\mbox{rng}}~ A$ denotes the closure of the range of $A$.\\
 $\dotplus$ (resp., $\oplus$) denotes the direct sum (resp., orthogonal direct sum) between two subspaces and $\ominus$ denotes the orthogonal complement.\\
A function $g:\mathbb{R}\to\mathfrak{X}$ is said to be integrable if it is Bochner integrable and square integrable if it satisfies the following condition 
$$\int_{-\infty}^\infty ||g(t)||^2dt<\infty.$$
The Fourier transformation of a square integrable function $g:\mathbb{R}\to\mathfrak{X}$ is denoted by $\hat{g}$ and is defined as
$$\hat{g}(t)=\int_{-\infty}^\infty e^{-i s t}g(s)ds.$$

\section{Operator valued  kernels and vector valued reproducing kernel Hilbert spaces (RKHS)}
\label{Se 2}
This section briefly recalls a few basic facts regarding vector valued RKHS. Since our goal is to work with RKHS of vector valued entire functions, we present all the results by assuming $\mathcal{H}$ to be the RKHS of $\mathfrak{X}$-valued entire functions. A detailed and general study can be found in \cite{Paulsen}.\\
Let $\mathcal{H}$ be a Hilbert space of $\mathfrak{X}$-valued entire functions. Then we call $\mathcal{H}$ a reproducing kernel Hilbert space if there exists a $B(\mathfrak{X})$-valued function $K_\xi(z)$ on $\mathbb{C}\times\mathbb{C}$, which satisfies the following two conditions: 
\begin{enumerate}
\item $K_\xi u\in\mathcal{H}$ for all $\xi\in\mathbb{C}$ and $u\in \mathfrak{X}$.
\item $\langle f,K_\xi u\rangle_\mathcal{H}=\langle f(\xi),u\rangle_\mathfrak{X}$ for all $f\in\mathcal{H}$, $\xi\in\mathbb{C}$ and $u\in\mathfrak{X}$.
\end{enumerate}
The $B(\mathfrak{X})$-valued function $K_\xi(z)$ is known as reproducing kernel (RK) for $\mathcal{H}$. 
Equivalently, $\mathcal{H}$ is an RKHS if for all $\xi\in\mathbb{C}$, the point evaluations 
$$\delta_\xi :\mathcal{H}\to \mathfrak{X},\hspace{.5cm} f\mapsto f(\xi)$$
are bounded. The function $L_\xi(z)=\delta_z\delta^*_\xi$  satisfies the two conditions of a reproducing kernel. For an RKHS, the reproducing kernel is unique. The supplementary calculation after assuming two reproducing kernels $K_\xi(z)$ and $G_\xi(z)$ for $\mathcal{H}$
\begin{align*}
||K_\xi u-G_\xi u||^2_\mathcal{H} & =\langle K_\xi u-G_\xi u,K_\xi u\rangle_\mathcal{H}-\langle K_\xi u-G_\xi u,G_\xi u\rangle_\mathcal{H}\\
& =\langle K_\xi(\xi)u-G_\xi(\xi)u,u\rangle_\mathfrak{X}-\langle K_\xi(\xi)u-G_\xi(\xi)u,u\rangle_\mathfrak{X}\\
& =0
\end{align*}
serves to verify our last statement. Thus the $B(\mathfrak{X})$-valued function $K_\xi(z)=\delta_z\delta^*_\xi$ is the reproducing kernel for $\mathcal{H}$. Clearly, for $\xi,z\in\mathbb{C}$, $K_\xi(z)^*=K_z(\xi)$. In an RKHS $\mathcal{H}$ norm convergence of a sequence of functions  implies point wise convergence, i.e. if $\{g_n\}\subseteq \mathcal{H}$ converges to $g\in\mathcal{H}$ in the norm, then for every $z\in\mathbb{C}$, $g_n(z)\to g(z)$. The linear span of the collection $\{K_\xi u:\xi\in\mathbb{C},u\in\mathfrak{X}\}$ is dense in $\mathcal{H}$. If there exist countable sets $\{\xi_1,\xi_2,\ldots\}\subset \mathbb{C}$ and $\{u_1,u_2,\ldots\}\subset\mathfrak{X}$ such that $\{K_{\xi_n}u_n:n\in\mathbb{N}\}$ is an orthogonal basis of $\mathcal{H}$, then we shall say that the RKHS has the Kramer sampling property which will appear in section~\ref{Se 9}.\\
 The reproducing kernel $K_\xi(z)$ is positive in the sense that, for every choice of $n\in\mathbb{N}$, $\xi_1,\xi_2,\ldots,\xi_n\in\mathbb{C}$ and $u_1,u_2,\ldots,u_n\in\mathfrak{X}$ the following is true
 $$\sum_{l,m=1}^n \bigg\langle K_{\xi_m}(\xi_l)u_m,u_l \bigg\rangle_\mathfrak{X}=\left\lVert\sum_{l=1}^n\delta_{\xi_l}^*(u_l)\right\rVert_\mathcal{H}^2\geq 0.$$ 
Clearly, for every $\xi\in\mathbb{C}$, $K_\xi(\xi)\succeq 0$. The set $\{\xi\in\mathbb{C}:K_\xi(\xi)~\mbox{is invertible}\}$ is an open subset of $\mathbb{C}$. The proof of the last assertion is the same as Lemma $2.4$ in \cite{JFA}.
A $B(\mathfrak{X})$-valued function $K_\xi(z)$ on $\mathbb{C}\times\mathbb{C}$ is called a positive kernel if it is positive in the sense as mentioned above. The operator version of Moore's theorem (Theorem 6.12, \cite{Paulsen}) ensures that corresponding to every positive kernel, there exists a unique RKHS $\mathcal{H}$. The next lemma gives a criterion to construct RKHS of entire functions.
\begin{lemma}
Let $\mathcal{H}$ be an RKHS of $\mathfrak{X}$-valued functions defined on $\mathbb{C}$ with RK $K_\xi(z)$. Then $\mathcal{H}$ is an RKHS of $\mathfrak{X}$-valued entire functions if and only if $K_\xi(z)$ is an entire function in $z$ for all $\xi\in\mathbb{C}$ and $||K_\xi(\xi)||$ is bounded on every compact subset of $\mathbb{C}$.
\end{lemma}
\begin{proof}
The proof is similar to Lemma $5.6$ in \cite{ArD08}.
\end{proof}
Next, we give an example of an RKHS of $\mathfrak{X}$-valued entire functions, which is a vector version of the Paley-Wiener space.

\begin{defn}
An entire function $g:\mathbb{C}\to\mathfrak{X}$ is said to be of exponential type at most $a$ if for each $\epsilon >0$ there exists a constant $L$, independent of $\xi$ such that
$$||g(\xi)||_\mathfrak{X}\leq L\hspace{.1cm} e^{|\xi |(a+\epsilon)} .$$
\end{defn}
If $g:\mathbb{R}\to\mathfrak{X}$ is a square integrable function, vanishes outside the compact interval $[-a,a]$, for some $a>0$, then $\hat{g}$ and $g$ satisfy the Plancherel's theorem and $\hat{g}$ can be extended as an entire function
$$\hat{g}(\xi)=\int_{-\infty}^\infty e^{-i\xi t}g(t)dt$$
 which is of exponential type at most $a$. These follow from Theorem $1.8.2$ and Theorem $1.8.3$ of the book \cite{Arendt}.
\begin{examp}[Paley-Wiener spaces of vector valued entire functions]
\label{Paley}
For $a>0$, the set of $\mathfrak{X}$-valued entire functions
 \begin{equation}
 PW_a=\{\hat{g}: g \hspace{.1cm}\mbox{is square integrable and vanishes outside the interval}\hspace{.1cm}[-a,a]\}
 \end{equation}
 is a Hilbert space with respect to the inner product
 \begin{equation}
 \langle \hat{g},\hat{h}\rangle_{PW_a}=\int_{-\infty}^\infty \langle \hat{g}(t),\hat{h}(t)\rangle_\mathfrak{X}dt .
 \end{equation}
 Also $PW_a$ is an RKHS with the reproducing kernel 
 $$K_\xi^a(z)=\frac{\sin(z-\overline{\xi})a}{\pi(z-\overline{\xi})}I_\mathfrak{X} ,$$
 where $I_\mathfrak{X}$ is the identity operator on $\mathfrak{X}$.
 Since for every $u\in\mathfrak{X}$ and $\xi\in\mathbb{C}$
 $$K_\xi^a(z) u=\int_{-\infty}^\infty e^{-i z t}Q_\xi^a(t) dt,$$
 where $Q_\xi^a(t)$ is a square integrable function defined as
 \begin{equation}
Q_\xi^a(t) := \left\{
\begin{array}{ll}
\frac{1}{2\pi}e^{i \overline{\xi} t}u  & \mbox{if } |t| \leq a \\
0 & \mbox{otherwise } 
\end{array} \right.
\end{equation}
and for $\hat{g}\in PW_a$,
\begin{align}
\langle \hat{g},K_\xi^a u \rangle_{PW_a} & = \langle \hat{g},\widehat{Q_\xi^a u} \rangle_{PW_a}\nonumber\hspace{.5cm}[~\because K_\xi^a u=\widehat{Q_\xi^a u}~] \\
& = 2\pi\langle g,Q_\xi^a u \rangle_{L^2}\nonumber\hspace{.5cm} [~\mbox{by Plancherel's theorem}~]\\
& = 2\pi \int_{-a}^a \langle g(t),\frac{1}{2\pi}e^{i \overline{\xi} t}u \rangle_{\mathfrak{X}}dt\nonumber\\
& = \left\langle \int_{-a}^a g(t)e^{-i \xi t}dt,u \right\rangle_{\mathfrak{X}}\nonumber\\
& = \langle \hat{g}(\xi),u\rangle_\mathfrak{X} . \nonumber
\end{align}
\end{examp}
\section{Spaces of vector valued holomorphic functions}
\label{Se 3}
In this section, we recall some crucial spaces of vector valued holomorphic functions. As most of the results are well known, we mention them here without proof. Details about these spaces can be found in \cite{Rosenblum}, \cite{Nagy}. As we have considered earlier, $\mathfrak{X}$ is a complex seperable Hilbert space, and $B(\mathfrak{X})$ is the algebra of bounded linear operators on $\mathfrak{X}$. We denote
$$L^2_\mathfrak{X}(\mathbb{R}):=\left\{f:\mathbb{R}\to\mathfrak{X}~|~f\hspace{.1cm}\mbox{is weakly measurable and}\hspace{.1cm}\int_{-\infty}^\infty ||f(t)||_\mathfrak{X}^2~ dt<\infty \right\},$$
$$H^2_\mathfrak{X}(\mathbb{C}_+):=\left\{f:\mathbb{C_+}\to\mathfrak{X}~|~f\hspace{.1cm}\mbox{is holomorphic and}~\mbox{sup}_{y>0}\int_{-\infty}^\infty ||f(x+iy)||^2_\mathfrak{X}~dx<\infty\right\}$$
and
$$H^\infty_{B(\mathfrak{X})}(\mathbb{C}_+):=\left\{f:\mathbb{C_+}\to B(\mathfrak{X})~|~f\hspace{.1cm}\mbox{is bounded and holomorphic}\right\}.$$
It is known that :-
\begin{enumerate}
 \item $L^2_\mathfrak{X}(\mathbb{R})$ is a Hilbert space with respect to the inner product
$$\langle f,g\rangle_{L^2}=\int_{-\infty}^\infty\langle f(t),g(t)\rangle_\mathfrak{X}dt$$
for all $f,g\in L^2_\mathfrak{X}(\mathbb{R})$.
\item The Hardy space over the upper half-plane $H^2_\mathfrak{X}(\mathbb{C}_+)$ is a Hilbert space with respect to the inner product
$$\langle f,g\rangle_{H^2}=\int_{-\infty}^\infty\langle f_0(x),g_0(x)\rangle_\mathfrak{X}dx$$
where $f_0,g_0\in L^2_\mathfrak{X}(\mathbb{R})$ are the boundary functions of $f$ and $g$ respectively, which are mentioned in the next theorem. 
\item $H^\infty_{B(\mathfrak{X})}(\mathbb{C}_+)$ is a Banach space with the norm
$$||F||_\infty=\mbox{sup}_{y>0}||F(x+iy)||_{B(\mathfrak{X})}$$
for all $F\in H^\infty_{B(\mathfrak{X})}(\mathbb{C}_+)$.
\end{enumerate}
We also denote 
$$\mathcal{S}=\left\{F\in H^\infty_{B(\mathfrak{X})}(\mathbb{C}_+): ||F(z)||\leq 1\hspace{.1cm}\mbox{for all}\hspace{.1cm}z\in\mathbb{C}_+\right\}.$$
The following two theorems give the boundary behaviour of the functions, respectively, in $H^2_\mathfrak{X}(\mathbb{C}_+)$ and $H^\infty_{B(\mathfrak{X})}(\mathbb{C}_+)$.
\begin{thm}
\label{Fatou's theorem}
Let $f\in H^2_\mathfrak{X}(\mathbb{C}_+)$, then there exists a (unique) nontangential boundary limit $f_0\in L^2_\mathfrak{X}(\mathbb{R})$ such that 
$$f_0(x)=\lim_{y\downarrow 0}f(x+iy)$$
 pointwise a.e. on $\mathbb{R}$. Also, $f_0$ satisfies the following identities
\begin{equation}
f(z)=\frac{1}{2\pi i}\int_{-\infty}^\infty \frac{f_0(t)}{t-z}dt\hspace{.5cm},y>0\label{1}
\end{equation}
and
\begin{equation}
\frac{1}{2\pi i}\int_{-\infty}^\infty \frac{f_0(t)}{t-z}dt=0\hspace{.5cm},y<0\label{2}
\end{equation}
where $z=x+iy$.\\
Conversely, every $f_0\in L^2_\mathfrak{X}(\mathbb{R})$ satisfying (\ref{1}) and (\ref{2}) gives the corresponding function $f\in H^2_\mathfrak{X}(\mathbb{C}_+)$ such that $f_0$ is the boundary function of $f$.
\end{thm}
The identity in (\ref{1}) is known as the Cauchy integral formula.
\begin{thm}
\label{Boundary value}
If $F\in H^\infty_{B(\mathfrak{X})}(\mathbb{C}_+)$, then for a.e. $x\in \mathbb{R}$ there exists $F_0(x)\in B(\mathfrak{X})$ such that for all $u\in \mathfrak{X}$
$$F(x+iy)u\to F_0(x)u\hspace{.5cm}\mbox{as}\hspace{.1cm}y\downarrow 0$$
and 
$$||F_0(x)||=\lim_{y\downarrow 0}||F(x+iy)||.$$
\end{thm}
We denote $\mathcal{S}^{in}$ (resp., $\mathcal{S}^{in}_*$) as the collection of all functions $F\in \mathcal{S}$ such that the corresponding boundary function $F(x)\in B(\mathfrak{X})$ is an isometry (resp., co-isometry) for a.e. $x\in\mathbb{R}$. It is easy to observe that a $B(\mathfrak{X})$-valued holomorphic function $F(z)$ on $\mathbb{C}_+$ belongs to $\mathcal{S}^{in}$ (resp., $\mathcal{S}^{in}_*$) if and only if $I_\mathfrak{X}-F(z)^*F(z)\succeq 0$ (resp., $I_\mathfrak{X}-F(z)F(z)^*\succeq 0$) for all $z\in\mathbb{C}_+$ with equality a.e. on $\mathbb{R}$. The operator valued functions $F\in \mathcal{S}^{in}$ (resp., $\mathcal{S}^{in}_*$) are called inner (resp., $*$-inner) functions. Functions belong to both $\mathcal{S}^{in}$ and $\mathcal{S}^{in}_*$ are called inner from both sides. \\
It is known that the characteristic function of a contraction operator $A\in B(\mathfrak{X})$,
$$C_A(z)=\left.\left[-A+z(I-AA^*)^\frac{1}{2}(I-zA^*)^{-1}(I-A^*A)^\frac{1}{2}\right]\right|{\overline{\mbox{rng}(I-A^*A)^\frac{1}{2}}}$$
 is an inner (resp., $*$-inner) function on the disc if and only if $A^{*n}\to 0$ (resp., $A^{n}\to 0$) strongly as $n\to\infty$ (see \cite{Nagy}). 
 Henry Helson studied in \cite{Helson}, the inner functions $F(z)$ from both sides, which are norm differentiable on the real line and satisfy the following differential equation
$$F'(x)=i~V(x)~F(x),$$
where $V(x)$ is $B(\mathfrak{X})$-valued norm continuous function, $V(x)\succeq 0$ and self adjoint for all $x\in\mathbb{R}$. Also, Stephen L. Campbell studied $B(\mathfrak{X})$-valued inner functions (see \cite{Campbell}), which are analytic on the closed unit disc. If $A\in B(\mathfrak{X})$ with the spectral radius $r(A)<1$, $||A||\leq 1$ and $AA^*\neq I_\mathfrak{X}$, then the corresponding Potapov inner function
$$V_A(z)=-A^*+z(I_\mathfrak{X}-A^*A)^\frac{1}{2}(I_\mathfrak{X}-zA)^{-1}(I_\mathfrak{X}-AA^*)^\frac{1}{2}$$
 is also analytic on the closed unit disc. We can consider inner functions on the disc as an inner function on the upper half-plane by using the conformal map $C(z)=\frac{z-i}{z+i}$ between the upper half-plane and the disc.\\
We denote $H^2_\mathfrak{X}(\mathbb{R})$ as the collection of all nontangential boundary limits of elements of $H^2_\mathfrak{X}(\mathbb{C}_+)$. Thus we can consider $H^2_\mathfrak{X}(\mathbb{C}_+)$ as a closed subspace of $L^2_\mathfrak{X}(\mathbb{R})$ in terms of $H^2_\mathfrak{X}(\mathbb{R})$. If we consider Hardy space over the lower half-plane, i.e., $H^2_\mathfrak{X}(\mathbb{C}_-)$, a similar result as Theorem \ref{Fatou's theorem} will also hold. The identity in (\ref{2}) implies that the orthogonal complement of $H^2_\mathfrak{X}(\mathbb{C}_+)$ can be identified with $H^2_\mathfrak{X}(\mathbb{C}_-)$.
\begin{thm}
The Hardy space over the upper half-plane $H^2_\mathfrak{X}(\mathbb{C}_+)$ and over the lower half-plane $H^2_\mathfrak{X}(\mathbb{C}_-)$ are RKHS of $\mathfrak{X}$-valued holomorphic functions on $\mathbb{C}_+$ and $\mathbb{C}_-$ respectively. The corresponding reproducing kernels are
\begin{equation}
K_\xi(z)=\frac{I_\mathfrak{X}}{\rho_\xi(z)}\hspace{.5cm}\xi,z\in\mathbb{C}_+
\end{equation}
and 
\begin{equation}
K_\alpha^{(-)}(\lambda)=-\frac{I_\mathfrak{X}}{\rho_\alpha(\lambda)}\hspace{.5cm}\alpha,\lambda\in\mathbb{C}_-.
\end{equation}
\end{thm}
\begin{proof}
A simple calculation shows that for every $u\in\mathfrak{X}$ and $\xi\in\mathbb{C}_+$, $K_\xi u\in H^2_\mathfrak{X}(\mathbb{C}_+)$. Also, from the Cauchy integral formula, the reproducing property 
\begin{align*}
\langle g,K_\xi u\rangle_{H_2}& =\int_{-\infty}^\infty\left\langle g_0(x),\frac{1}{-2\pi i(x-\overline{\xi})}u\right\rangle_\mathfrak{X}dx \\
& =\left\langle \frac{1}{2\pi i}\int_{-\infty}^\infty \frac{g_0(x)}{x-\xi}dx,u\right\rangle_\mathfrak{X}\\
& =\langle g(\xi),u\rangle_\mathfrak{X} \nonumber
\end{align*}
holds for any $g\in H^2_\mathfrak{X}(\mathbb{C}_+)$, $\xi\in\mathbb{C}_+$ and $u\in\mathfrak{X}$. In the similar way the case for $H^2_\mathfrak{X}(\mathbb{C}_-)$ can also be proved.
\end{proof}
Suppose $F\in H^\infty_{B(\mathfrak{X})}(\mathbb{C}_+)$ with $||F(z)||\leq M$ for all $z\in\mathbb{C}_+$, then we can define an operator $M_F:H^2_\mathfrak{X}(\mathbb{C}_+)\to H^2_\mathfrak{X}(\mathbb{C}_+)$ defined by 
\begin{equation}
(M_Fg)(z)=F(z)g(z)\hspace{.5cm}\mbox{for}\hspace{.1cm}g\in H^2_\mathfrak{X}(\mathbb{C}_+).\label{3}
\end{equation}
The following evaluation shows that $M_F$ is well defined.
For $g\in H^2_\mathfrak{X}(\mathbb{C}_+)$
\begin{align*}
\int_{-\infty}^\infty ||F(x+iy)g(x+iy)||_\mathfrak{X}^2 dx & \leq \int_{-\infty}^\infty ||F(x+iy)||^2||g(x+iy)||_\mathfrak{X}^2dx\\ 
& \leq M^2 \int_{-\infty}^\infty ||g(x+iy)||_\mathfrak{X}^2dx\\
& \leq M^2 ||g||^2_{H^2_\mathfrak{X}(\mathbb{C}_+)}.
\end{align*}
\begin{prop}
If $M_F$ is defined as (\ref{3}) for some $F\in H^\infty_{B(\mathfrak{X})}(\mathbb{C}_+)$, then the following implications are true:
\begin{enumerate}
\item $M_F$ is a bounded operator on $H^2_\mathfrak{X}(\mathbb{C}_+)$ with $||M_F||=||F||_\infty$.
\item $M_F^*\frac{u}{\rho_\xi}=F(\xi)^*\frac{u}{\rho_\xi}$ for all $u\in\mathfrak{X}$ and $\xi\in\mathbb{C}_+$.
\end{enumerate}
Moreover, if $F\in \mathcal{S}$ then $M_F$ is a contractive operator.
\end{prop}
\begin{proof}
Boundedness of $M_F$ and $||M_F||\leq ||F||_\infty$ follows from the preceding calculations. Also a simple use of the Cauchy integral formula shows that for all $u,v\in\mathfrak{X}$ and $\xi\in\mathbb{C}_+$ the following inequality is true
$$|\langle F(\xi)u,v\rangle_\mathfrak{X}|\leq ||M_F||\hspace{.1cm}||u||\hspace{.1cm}||v||.$$
This implies for all $\xi\in\mathbb{C}_+$, $||F(\xi)||\leq ||M_F||$. Therefore,  $||M_F||=||F||_\infty$.\\
Since the RK for $H^2_\mathfrak{X}(\mathbb{C}_+)$ is $K_\xi(z)=\frac{I_\mathfrak{X}}{\rho_\xi(z)}$. To show $(2)$ it is sufficient to show that for all $u,v\in\mathfrak{X}$ and $\xi,z\in\mathbb{C}_+$
$$\langle M_F \frac{u}{\rho_\xi},\frac{v}{\rho_z}\rangle_{H^2}=\langle \frac{u}{\rho_\xi},F(z)^* \frac{v}{\rho_z}\rangle_{H^2},$$
which can be shown by a simple calculation using the Cauchy integral formula. If $F\in \mathcal{S}$ then $M_F$ is a contractive operator follows from $(1)$. 
\end{proof}
In view of the preceding proposition, for any $F\in \mathcal{S}$ and any $n\in\mathbb{N}$
\begin{equation}
\sum_{l,m=1}^n\left\langle \frac{I_\mathfrak{X}-F(z_l)F(z_m)^*}{\rho_{z_m}(z_l)}u_m,u_l\right\rangle_\mathfrak{X}\geq 0
\end{equation}
for every choice of $u_1,u_2,\ldots,u_n\in\mathfrak{X}$ and $z_1,z_2,\ldots,z_n\in\mathbb{C}_+$. Thus the $B(\mathfrak{X})$-valued function 
$$\Gamma_\xi(z)=\frac{I_\mathfrak{X}-F(z)F(\xi)^*}{\rho_\xi(z)}$$
is a positive kernel on $\mathbb{C}_+\times\mathbb{C}_+$, and we denote the corresponding unique RKHS of $\mathfrak{X}$-valued holomorphic functions on $\mathbb{C}_+$ as $\mathcal{H}(F)$.
\section{Characterization and extension of the RKHS $\mathcal{H}(F)$}
\label{Se 4}
In this section, we recall an analogues characterization of the space $\mathcal{H}(F)$, which has been mentioned in \cite{Branges 4}, and extend $\mathcal{H}(F)$ as an RKHS of $\mathfrak{X}$-valued holomorphic functions on a domain possibly larger than $\mathbb{C}_+$. We denote $\mathcal{P}$ as the orthogonal projection of $L^2_\mathfrak{X}(\mathbb{R})$ onto $H^2_\mathfrak{X}(\mathbb{R})$ and $\mathcal{Q}=I_{L^2_\mathfrak{X}(\mathbb{R})}-\mathcal{P}$. For $F\in \mathcal{S}$ and $f\in H^2_\mathfrak{X}(\mathbb{C}_+)$, we denote
$$\nu(f)=\mbox{sup}\left\{||f+M_F(g)||^2_{H^2}-||g||^2_{H^2}:g\in H^2_\mathfrak{X}(\mathbb{C}_+)\right\}.$$
\begin{thm}
\label{Charac}
For $F\in \mathcal{S}$
$$\mathcal{H}(F)=\left\{f\in H^2_\mathfrak{X}(\mathbb{C}_+):\nu(f)<\infty\right\}\hspace{.1cm}\mbox{and}\hspace{.1cm}||f||^2_{\mathcal{H}(F)}=\nu(f).$$
Moreover, if $F\in \mathcal{S}^{in}$, then
$$\mathcal{H}(F)=H^2_\mathfrak{X}(\mathbb{C}_+)\ominus \mbox{rng}M_F\hspace{.1cm}\mbox{and}\hspace{.1cm} ||f||^2_{\mathcal{H}(F)}=||f||^2_{H^2}.$$
\end{thm}
\begin{proof}
For proof, see Theorem $2.2$ and its second corollary in \cite{Dym89}.
\end{proof}
In view of Theorem \ref{Boundary value}, for every $F\in \mathcal{S}$ and $b\in L^2_\mathfrak{X}(\mathbb{R})$ the function defined by
$$x\mapsto F_0(x)b(x)\hspace{.5cm}\mbox{for}\hspace{.1cm} x\in\mathbb{R}$$
belongs to $L^2_\mathfrak{X}(\mathbb{R})$, where $F_0$ is the nontangential boundary limit of $F$. We denote this function as $Fb$.
\begin{lemma}
Let $F\in \mathcal{S}$, $b\in L^2_\mathfrak{X}(\mathbb{R})\ominus H^2_\mathfrak{X}(\mathbb{R})$ and $f\in H^2_\mathfrak{X}(\mathbb{C}_+)$ be such that $f_0=-\mathcal{P}Fb$ is the corresponding nontangential boundary limit. Then for all $g\in H^2_\mathfrak{X}(\mathbb{C}_+)$
$$ ||f+M_F(g)||^2-||g||^2\leq ||b||^2-||\mathcal{Q}Fb||^2,$$
with equality if $F\in \mathcal{S}^{in}$.
\end{lemma}
\begin{proof}
Since $f$ and $M_F(g)$ both belong to $H^2_\mathfrak{X}(\mathbb{C}_+)$,
$$||f+M_F(g)||_{H^2}=||f_0+Fg_0||_{L^2},$$
where $g_0$ is the nontangential boundary limit of $g$. Now,
$$||-\mathcal{P}Fb+Fg_0||^2_{L^2}= ||\mathcal{Q}Fb+F(g_0-b)||^2_{L^2}=-||\mathcal{Q}Fb||^2_{L^2}+||F(g_0-b)||^2_{L^2}$$
as $\langle \mathcal{Q}Fb,Fg_0\rangle_{L^2}=0$ and $\langle \mathcal{Q}Fb,Fb\rangle_{L^2}=||\mathcal{Q}Fb||^2_{L^2}$. Thus the lemma follows from the following fact that 
$$||F(g_0-b)||^2_{L^2}\leq ~\mbox{or} ~ = ||g_0||^2_{L^2}+||b||^2_{L^2}$$
according as $F$ belongs to $ \mathcal{S}$  or $ \mathcal{S}^{in}$.
\end{proof}
The next theorem is motivated from \cite{Dym89} (Theorem $2.3$), where we replaced the matrix valued kernels with operator valued kernels. Since these results are crucial for the rest of this article, we include the proof with suitable modifications.
\begin{thm}
\label{Inequality}
If $F\in \mathcal{S}$, then for every choice of $n\in\mathbb{N}$, $z,z_1,z_2,\ldots,z_n\in\mathbb{C}_+$ and $u,u_1,u_2,\ldots,u_n\in\mathfrak{X}$ the following implications hold:
\begin{enumerate}
\item $(R_zF)u\in\mathcal{H}(F)$.
\item $||\sum_{l=1}^n(R_{z_l}F)u_l||^2_{\mathcal{H}(F)}\leq 4\pi^2\sum_{l,m=1}^n\left\langle\frac{I_\mathfrak{X}-F(z_m)^*F(z_l)}{\rho_{z_m}(z_l)}u_l,u_m\right\rangle_\mathfrak{X}$ with equality if $F\in \mathcal{S}^{in}$.
\item $\mathcal{H}(F)$ is invariant under $R_z$ for all $z\in\mathbb{C}_+$.
\item $R_z$ is a bounded operator on $\mathcal{H}(F)$ for all $z\in \mathbb{C}_+$ and for all $g\in \mathcal{H}(F)$ $R_z$ satisfy the following assertion
\begin{equation}
||R_zg||^2\leq \frac{\mbox{Im}(\langle R_zg,g\rangle)-\pi ||g(z)||^2}{\mbox{Im}(z)}.\label{7}
\end{equation}
\end{enumerate}
\end{thm}
\begin{proof}
Consider $$b=\sum_{l=1}^n\frac{u_l}{\xi-z_l}.$$
Clearly, $b\in H^2_\mathfrak{X}(\mathbb{C}_-)$. Thus it's nontangential boundary function $b$ (say) belongs to $L^2_\mathfrak{X}(\mathbb{R})\ominus H^2_\mathfrak{X}(\mathbb{C}_+)$. Now 
$$(R_zF)(\xi)u=\frac{F(\xi)u-F(z)u}{\xi-z}$$
is analytic in $\mathbb{C}_+$, and the nontangential boundary function $(R_zF)(t)u$ belongs to $H^2_\mathfrak{X}(\mathbb{R})$ as
$$\frac{1}{2\pi i}\int_{-\infty}^\infty \frac{(R_zF)(t)u}{t-\alpha}dt=0$$ 
for all $\alpha\in\mathbb{C}_-$.
Since $F(t)b(t)\in L^2_\mathfrak{X}(\mathbb{R})$ and
$$F(t)b(t)=\left(\sum_{l=1}^n \frac{F(t)u_l-F(z_l)u_l}{t-z_l}\right)+\left(\sum_{l=1}^n \frac{F(z_l)u_l}{t-z_l}\right),$$ then
$$\sum_{l=1}^n (R_{z_l}F)u_l=\mathcal{P}Fb$$
and
$$\sum_{l=1}^n \frac{F(z_l)u_l}{t-z_l}=\mathcal{Q}Fb.$$
Now applying the preceding lemma, we have for all $g\in H^2_\mathfrak{X}(\mathbb{C}_+)$
\begin{equation}
||\sum_{l=1}^n (R_{z_l}F)u_l+M_F(g)||^2-||g||^2\leq ||b||^2-||\mathcal{Q}Fb||^2.\label{4}
\end{equation}
Thus $(1)$ follows from (\ref{4}) in association with Theorem \ref{Charac}. \\
A simple calculation by using the Cauchy integral formula gives
\begin{equation}
||b||^2-||\mathcal{Q}Fb||^2= 4\pi^2\sum_{l,m=1}^n\left\langle\frac{I_\mathfrak{X}-F(z_m)^*F(z_l)}{\rho_{z_m}(z_l)}u_l,u_m\right\rangle_\mathfrak{X}.\label{5}
\end{equation}
Thus $(2)$ follows from (\ref{4}) with $g\equiv 0$ and (\ref{5}).\\
Since for every $z\in\mathbb{C}_+$, $R_z$ is linear, and $\mathcal{H}(F)$ is an RKHS, to show $(3)$ and $(4)$ it is sufficient to show that they are valid for $\Gamma_\alpha u$ for every choice of $u\in\mathfrak{X}$ and $\alpha\in\mathbb{C}_+$. Since
\begin{equation}
R_z\Gamma_\alpha(\xi)u=\frac{2\pi i}{\rho_\alpha(z)}\Gamma_\alpha(\xi)u-\frac{1}{\rho_\alpha(z)}(R_z F)(\xi)F(\alpha)^*u\label{6}
\end{equation}
and both the terms in RHS of (\ref{6}) belongs to $\mathcal{H}(F)$. Thus $R_z\Gamma_\alpha u\in \mathcal{H}(F)$ and $\Gamma_\alpha u$ satisfies (\ref{7}).
\end{proof}
Suppose $F\in \mathcal{S}$ is inner from both sides and consider $\mathfrak{A}_F^-=\{z\in\mathbb{C}_-:F(\overline{z})\hspace{.1cm}\mbox{is invertible}\}$. Since $F$ is invertible at one point implies it is invertible in a certain open neighbourhood of that point, $\mathfrak{A}_F^-$ is an open subset of $\mathbb{C}_-$. Now we can extend $F$ as a $B(\mathfrak{X})$-valued holomorphic function on $\mathfrak{A}_F^-$ by the following formula
\begin{equation}
F(z)=\{F(\overline{z})^*\}^{-1}\hspace{.5cm}\mbox{for}\hspace{.1cm}z\in \mathfrak{A}_F^-.\label{8}
\end{equation}
Also, for almost every $x\in\mathbb{R}$
$$F(x)=\lim_{y\downarrow 0}F(x+iy)=\lim_{y\downarrow 0}F(x-iy).$$
Now, for a function $F\in \mathcal{S}$ and inner from both sides, we denote $\mathfrak{F}$ as the extension of $F$ and $\mathfrak{A}_\mathfrak{F}$ (containing $\mathbb{C}_+\cup \mathfrak{A}_F^-$) as the domain of holomorphy of $\mathfrak{F}$.\\
Once we have the above extension $\mathfrak{F}$ of $F$, we can think $\mathcal{H}(F)$ as the RKHS of $\mathfrak{X}$-valued holomorphic functions on $\mathfrak{A}_\mathfrak{F}$, which we denote as $\mathcal{H}(\mathfrak{F})$. The following lemma provides more details of $\mathcal{H}(\mathfrak{F})$.
\begin{lemma}
\label{Extension}
If $F\in \mathcal{S}$ is inner from both sides and $\mathfrak{F}$ is the corresponding extension of $F$ as defined in (\ref{8}), the kernel function $K_\xi^{\mathfrak{F}}(z)$ defined by
\begin{equation}
     K_\xi^{\mathfrak{F}}(z):= \left\{
    \begin{array}{ll}
         \frac{I_\mathfrak{X}-\mathfrak{F}(z) \mathfrak{F}(\xi)^*}{\rho_\xi(z)}  & \mbox{if } z \neq \overline{\xi} \\
         \frac{\mathfrak{F}'(\overline{\xi})\mathfrak{F}(\xi)^*}{2\pi i} & \mbox{if } z = \overline{\xi}
    \end{array} \right.
\end{equation}
on $\mathfrak{A}_\mathfrak{F}\times\mathfrak{A}_\mathfrak{F}$ is positive.
\end{lemma}
\begin{proof}
To show that $K_\xi^{\mathfrak{F}}(z)$ is positive, we need to show that for every choice of $n\in\mathbb{N}$, $z_1,z_2,\ldots,z_n\in\mathfrak{A}_\mathfrak{F}$ and $u_1,u_2,\ldots,u_n\in\mathfrak{X}$
\begin{equation}
\sum_{l,m=1}^n\left\langle K_{z_m}^\mathfrak{F}(z_l)u_m,u_l\right\rangle_\mathfrak{X}\geq 0.\label{9}
\end{equation}
Here we only show the case where some points belong to $\mathbb{C}_+\cap\mathfrak{A}_\mathfrak{F}$ and others belong to $\mathbb{C}_-\cap\mathfrak{A}_\mathfrak{F}$. The remaining cases can be deduce from this. Without loss of generality we may assume that $z_1,z_2,\ldots,z_i\in \mathbb{C}_+\cap\mathfrak{A}_\mathfrak{F}$ and $\xi_1,\xi_2,\ldots,\xi_j\in \mathbb{C}_-\cap\mathfrak{A}_\mathfrak{F}$ with $i+j=n$ and $\xi_k=z_{i+k}$ for every $k=1,2,\ldots,j$. Also we assume that $v_k=u_{i+k}$ for all $k=1,2,\ldots,j$.\\
Since $\mathfrak{F}(\overline{\xi})\mathfrak{F}(\xi)^*=I_\mathfrak{X}$, for $z\neq\overline{\xi}$ the followings are true
\begin{align}
K_\xi^\mathfrak{F}(z) & = \frac{I_\mathfrak{X}-\mathfrak{F}(z)\mathfrak{F}(\xi)^*}{\rho_\xi(z)}\nonumber\\
& = \frac{\{\mathfrak{F}(\overline{\xi})-\mathfrak{F}(z)\}\mathfrak{F}(\xi)^*}{-2\pi i(z-\overline{\xi})}\nonumber\\
& = \frac{1}{2\pi i}(R_{\overline{\xi}}\mathfrak{F})(z)\mathfrak{F}(\xi)^*\label{10}
\end{align}
and
\begin{equation}
K_\xi^\mathfrak{F}(z)^*=\frac{1}{2\pi i}(R_{\overline{z}}\mathfrak{F})(\xi)\mathfrak{F}(z)^*.\label{11}
\end{equation}
Now we divide the sum in (\ref{9}) into four parts and simplify them with the help of (\ref{10}) and (\ref{11}).
The first part is
\begin{equation}
\sum_{l,m=1}^i\langle K_{z_l}^\mathfrak{F}(z_m)u_l,u_m\rangle_\mathfrak{X}.
\end{equation}
The second part is
\begin{align}
\sum_{t=1}^j\sum_{l=1}^i\langle K_{z_l}^\mathfrak{F}(\xi_t)u_l,v_t\rangle_\mathfrak{X} & = \sum_{t=1}^j\sum_{l=1}^i \langle  u_l,\frac{1}{2\pi i}(R_{\overline{\xi_t}}\mathfrak{F})(z_l)\mathfrak{F}(\xi_t)^* v_t\rangle_\mathfrak{X}\nonumber\\
& = \sum_{t=1}^j\sum_{l=1}^i\langle u_l,(R_{\overline{\xi_t}}\mathfrak{F})(z_l)x_t\rangle_\mathfrak{X},
\end{align}
where $x_t=\frac{1}{2\pi i}\mathfrak{F}(\xi_t)^*v_t$ for all $t=1,2,\ldots,j$.\\
The third part is
\begin{align}
\sum_{m=1}^j\sum_{s=1}^i\langle K_{\xi_m}^\mathfrak{F}(z_s)v_m,u_s\rangle_\mathfrak{X} & = \sum_{m=1}^j\sum_{s=1}^i\langle \frac{1}{2\pi i}(R_{\overline{\xi_m}}\mathfrak{F})(z_s)\mathfrak{F}(\xi_m)^* v_m,u_s\rangle_\mathfrak{X}\nonumber\\
& = \sum_{m=1}^j\sum_{s=1}^i\langle (R_{\overline{\xi_m}}\mathfrak{F})(z_s)x_m,u_s\rangle_\mathfrak{X}.
\end{align}
The fourth part is
\begin{align}
\sum_{t,m=1}^j\langle K_{\xi_m}^\mathfrak{F}(\xi_t)v_m,v_t\rangle_\mathfrak{X} & = \sum_{t,m=1}^j\langle \frac{I_\mathfrak{X}-\mathfrak{F}(\xi_t)\mathfrak{F}(\xi_m)^*}{\rho_{\xi_m}(\xi_t)}v_m,v_t\rangle_\mathfrak{X}\nonumber\\
& = \sum_{t,m=1}^j\langle \frac{\mathfrak{F}(\xi_t)\{\mathfrak{F}(\overline{\xi_t})^* \mathfrak{F}(\overline{\xi_m})-I_\mathfrak{X}\}\mathfrak{F}(\xi_m)^*}{-2\pi i(\xi_t-\overline{\xi_m})}v_m ,v_t \rangle_\mathfrak{X}\nonumber\\
& = \sum_{t,m=1}^j\langle\frac{ \{\mathfrak{F}(\overline{\xi_t})^* \mathfrak{F}(\overline{\xi_m})-I_\mathfrak{X}\}\mathfrak{F}(\xi_m)^*}{-2\pi i(\xi_t-\overline{\xi_m})}v_m,\mathfrak{F}(\xi_t)^* v_t\rangle_\mathfrak{X}\nonumber\\
& = 4\pi^2\sum_{t,m=1}^j\langle\frac{I_\mathfrak{X}-\mathfrak{F}(\overline{\xi_t})^*\mathfrak{F}(\overline{\xi_m})}{\rho_{\overline{\xi_t}}(\overline{\xi_m})}x_m,x_t\rangle_\mathfrak{X}.
\end{align}
In view of Theorem \ref{Inequality}, we have
\begin{equation}
\sum_{t,m=1}^j\langle K_{\xi_m}^\mathfrak{F}(\xi_t)v_m,v_t\rangle_\mathfrak{X}\geq ||\sum_{m=1}^j(R_{\overline{\xi_m}}\mathfrak{F})x_m||^2.
\end{equation}
Thus the above calculations ensure that
$$\sum_{l,m=1}^n\left\langle K_{z_m}^\mathfrak{F}(z_l)u_m,u_l\right\rangle_\mathfrak{X}\geq ||\sum_{l=1}^i K_{z_l}^\mathfrak{F}u_l+\sum_{m=1}^j(R_{\overline{\xi_m}}\mathfrak{F})x_m||^2\geq 0.$$
This completes the proof.
\end{proof}
Since for every $\xi\in\mathfrak{A}_\mathfrak{F}$, $K_\xi^\mathfrak{F}(\xi)$ is bounded, every $f\in\mathcal{H}(\mathfrak{F})$ are $\mathfrak{X}$-valued holomorphic functions on $\mathfrak{A}_\mathfrak{F}$.
\section{Construction of de Branges spaces based on pairs of Fredholm operator valued entire functions}
\label{Se 5}
This section consists of some of our main results. Our goal is to construct the de Branges operator $\mathfrak{E}=(E_-,E_+)$, likewise the de Branges matrix, as mentioned in the introduction. Here we shall see that the Fredholm operator valued holomorphic functions will play a significant role. In particular, it will be apparent that the components $E_\pm$ of the de Branges operator $\mathfrak{E}$ should be Fredholm operator valued entire functions with some additional properties to make the theory compatible with the existing theory of de Branges spaces based on de Branges matrices.
  Let us first see some basic results about Fredholm operators and Fredholm operator valued holomorphic functions. An operator $A \in B(\mathfrak{X})$ is said to be a Fredholm operator if it satisfies the following conditions:
\begin{enumerate}
\item $\dim(\ker(A))<\infty$.
\item $\mbox{rng(A)}$ is closed in $\mathfrak{X}$.
\item $\dim(\ker(A^*))<\infty$.
\end{enumerate}
We denote the collection of all Fredholm operators in $B(\mathfrak{X})$ as $\Phi(\mathfrak{X})$. For every $A\in \Phi(\mathfrak{X})$ the corresponding Fredholm index is defined by
$$\mbox{ind}(A)=\dim(\ker(A))-\dim(\ker(A^*)).$$
\begin{lemma}
If $A,B\in B(\mathfrak{X})$, then the following assertions are true
\begin{enumerate}
\item $A\in \Phi(\mathfrak{X}) \iff A^*\in \Phi(\mathfrak{X})$.
\item $A,B\in \Phi(\mathfrak{X})$ implies the composition $AB\in \Phi(\mathfrak{X})$.
\end{enumerate}
\end{lemma}
For more details about Fredholm operators, see \cite{Kato}. The next theorem is a particular form of analytic Fredholm theorem, which can be found in (Theorem $3.3$, \cite{Holden}, \cite{Gohberg1}) and references therein.
\begin{thm}
\label{Fredholm thm}
Suppose $\mathfrak{A}\subseteq \mathbb{C}$ be open and connected and $F:\mathfrak{A}\to B(\mathfrak{X})$ is analytic such that for all $z\in \mathfrak{A}$, $F(z)\in \Phi(\mathfrak{X})$. Then one of the two following assertions is always true
\begin{enumerate}
\item $F(z)^{-1}\not\in B(\mathfrak{X})$ for any $z\in \mathfrak{A}$.
\item $F(z)^{-1}\in \Phi(\mathfrak{X})$ for all $z\in \mathfrak{A}$, possibly except for a discrete set $D$. Also, the function $F(.)^{-1}$ is holomorphic on $\mathfrak{A}\setminus D$ and meromorphic on $\mathfrak{A}$.
\end{enumerate}
\end{thm}
Now we can start the process of constructing the new RKHS with operator valued RK. Let $E_+,E_-:\mathbb{C}\to  B(\mathfrak{X})$ be two entire functions such that $E_+(z), E_-(z)\in \Phi(\mathfrak{X})$ for all $z\in\mathbb{C}$. Also, assume that $E_+$ and $E_-$ both are invertible at least at one point and 
$$F:=E_+^{-1}E_-\in \mathcal{S}^{in}\cap\mathcal{S}^{in}_*.$$
Since $E_+$ and $E_-$ are invertible at least at one point, from the preceding theorem, we can find two discrete subsets, $D_1$ and $D_2$ of $\mathbb{C}$, such that $E_+$ is invertible for every $z\in \mathbb{C}\setminus D_1$ and $E_-$ is invertible for every $z\in \mathbb{C}\setminus D_2$. Also, $ F\in \mathcal{S}^{in}\cap\mathcal{S}^{in}_*$ implies that for every $z\in\mathbb{C}_+$
$$E_+(z)E_+(z)^*-E_-(z)E_-(z)^*\succeq 0$$
and for every $z\in\mathbb{R}$
$$E_+(z)E_+(\overline{z})^*-E_-(z)E_-(\overline{z})^*=0,$$
which can be extended to 
\begin{equation}
E_+(z)E_+(\overline{z})^*-E_-(z)E_-(\overline{z})^*=0\hspace{.5cm}\mbox{for every}\hspace{.1cm}z\in\mathbb{C}.\label{12}
\end{equation}
We call the pair of operator valued functions
\begin{equation}
\mathfrak{E}(z)=(E_-(z),E_+(z)) \hspace{.5cm}\mbox{for every}~z\in\mathbb{C}
\end{equation}
as de Branges operator. Now corresponding to the de Branges operator $\mathfrak{E}(z)$ we define the kernel 
\begin{equation}
K_\xi^{\mathfrak{E}}(z):= \left\{
    \begin{array}{ll}
         \frac{E_+(z)E_+(\xi)^*-E_-(z)E_-(\xi)^*}{\rho_\xi(z)}  & \mbox{if } z \neq \overline{\xi} \\
         \frac{E_+^{'} (\overline{\xi})E_+(\xi)^*- E_-^{'}(\overline{\xi})E_-(\xi)^*}{-2\pi i} & \mbox{if } z = \overline{\xi}
    \end{array} \right.\label{de Branges kernel}
\end{equation}
on $\mathbb{C}\times\mathbb{C}$.\\
Now we intend to show that the kernel defined in (\ref{de Branges kernel}) is positive on $\mathbb{C}\times\mathbb{C}$. Here we follow the process of extension as mentioned in the previous section. We denote the extended function of $F$ as $\mathfrak{F}$ and the extended domain as $\mathfrak{A}_\mathfrak{F}$. Observe that $\mathfrak{A}_\mathfrak{F}$ is dense in $\mathbb{C}$. Also, for all $\xi,z\in \mathfrak{A}_\mathfrak{F}$ 
\begin{equation}
K_\xi^{\mathfrak{E}}(z)=E_+(z)K_\xi^\mathfrak{F}(z)E_+(\xi)^*
\end{equation}
holds, which can be shown with the help of (\ref{12}). Since $K_\xi^\mathfrak{F}(z)$ is positive on $\mathfrak{A}_\mathfrak{F}\times \mathfrak{A}_\mathfrak{F}$ and $\mathfrak{A}_\mathfrak{F}$ is dense in $\mathbb{C}$, $K_\xi^{\mathfrak{E}}(z)$ is positive on $\mathbb{C}\times\mathbb{C}$. Thus, we can have a unique RKHS of $\mathfrak{X}$-valued entire functions corresponding to the positive definite kernel $K_\xi^{\mathfrak{E}}(z)$. We denote this RKHS as $\mathcal{B}(\mathfrak{E})$, and by observing the structure of the RK, we call it de Branges space. In particular, in a de Branges space $\mathcal{B}(\mathfrak{E})$, if for some $\beta\in\mathbb{C}_+$, $E_+(\beta)$ and $E_-(\overline{\beta})$ are self adjoint, then we denote the space as $\mathcal{B}_\beta(\mathfrak{E})$. A characterization of the space $\mathcal{B}_\beta(\mathfrak{E})$ can be found in Section~\ref{Se 8}.
\begin{rmk}
If $X\in B(\mathfrak{X})$ is any Fredholm operator and $XX^*=I$, then the pair of operators
$$\mathfrak{E}^X(z)=(E_-(z)X,E_+(z)X)\hspace{.5cm}\mbox{for every}~z\in\mathbb{C}$$ 
is again a de Branges operator and $\mathcal{B}(\mathfrak{E})=B(\mathfrak{E}^X)$.
\end{rmk}
\begin{rmk}
\label{index}
It is known that corresponding to the Fredholm operator valued holomorphic function $F(z)$ on a domain $\mathfrak{A}$ the index function defined by
$$z\mapsto \mbox{ind}F(z)\hspace{.5cm}\mbox{for all}~z\in\mathfrak{A}$$
is an integer valued continuous function. Thus corresponding to the de Branges operator $\mathfrak{E}(z)=(E_-(z),E_+(z))$, there exists a pair of integers. It is easy to observe that the pair of integers is always $(0,0)$ for de Branges operators.

\end{rmk}
 The following theorem provides the connection between $\mathcal{B}(\mathfrak{E})$  and the Hardy space.
\begin{thm}
If $\mathcal{B}(\mathfrak{E})$ is a de Branges space as defined above with respect to the de Branges operator $\mathfrak{E}(z)=(E_-(z),E_+(z))$ for every $z\in\mathbb{C}$. Then
\begin{equation}
\mathcal{B}(\mathfrak{E})=\{f:\mathbb{C}\to\mathfrak{X}~|~f~\mbox{is entire,}~E_+^{-1}f\in H^2_\mathfrak{X}(\mathbb{C}_+)~\mbox{and}~ E_-^{-1}f\in(H^2_\mathfrak{X}(\mathbb{C}_+))^\perp\}.
\end{equation}
Also with respect to the inner product
\begin{equation}
\langle f,g\rangle_{\mathcal{B}(\mathfrak{E})}=\int_{-\infty}^\infty \langle E_+^{-1}f(x),E_+^{-1}g(x)\rangle_\mathfrak{X}~dx,
\end{equation}
$\mathcal{B}(\mathfrak{E})$  is an RKHS, and the reproducing kernel is given by (\ref{de Branges kernel}).
\end{thm}
\begin{proof}
For proof, see Theorem $3.10$ in \cite{ArD18}.
\end{proof}

\section{Few examples}
\label{Se 6}
In this section, we present a few examples of the de Branges spaces which we have constructed in the last section.
\begin{examp}
 Consider $E_+(z)=e^{-iza}I_\mathfrak{X}$ and $E_-(z)=e^{iza}I_\mathfrak{X}$ for some $a>0$. Then, the corresponding de Branges space is actually the vector valued Paley-Wienner space as mentioned in example \ref{Paley}.
\end{examp}
  The next example is motivated by a Fredholm operator valued holomorphic function from \cite{Gohberg1} (chap. XI, sec. $2~ \& ~ 8$). Here we briefly mention this Fredholm operator valued holomorphic function.
   We denote $L_2([s,t])$ as the collection of all complex valued square integrable functions on $[s,t]$. Let us consider the boundary value problem
\begin{equation}
     \left\{
    \begin{array}{ll}
        g'(\mu)=P~g(\mu)+h(\mu),~~~s\leq \mu\leq t  \\
         Q_1~g(s)+Q_2~g(t)=u,
    \end{array} \right.\label{ex}
\end{equation}
where $h\in L_2^n([s,t])$, $u\in\mathbb{C}^n$ are given and $P, Q_1, Q_2$ are $n\times n$ matrices. The task is to find a solution in $G=\left(W_2^1([s,t])\right)^n$, where
\begin{equation}
W_2^1([s,t])=\{g\in L_2([s,t]): g~\mbox{is absolutely continuous}, g'\in L_2([s,t])\}
\end{equation}
is the Sobolev space of order one on $[s,t]$.
The operator form of (\ref{ex}) is
\begin{equation}
Ag:=
     \left[
    \begin{array}{ll}
        g'-P~g  \\
         Q_1~g(s)+Q_2~g(t)
    \end{array} \right]=
    \left[
    \begin{array}{ll}
        h  \\
         u
    \end{array} \right].
\end{equation}
In \cite{Gohberg1}, it has been proved that $A:G\to L_2^n([s,t])\oplus \mathbb{C}^n$ is a Fredholm operator of index zero.\\
Suppose $\mathcal{G}\subseteq \mathbb{C}$ is open and connected. Let $P(z)$, $Q_1(z)$ and $Q_2(z)$ are $n\times n$ matrix valued holomorphic functions on $\mathcal{G}$. The following boundary value problem gives rise to a Fredholm operator valued holomorphic function
\begin{equation}
\label{BP}
     \left\{
    \begin{array}{ll}
        g'(\mu)=P(z)~g(\mu)+h(\mu),~~~s\leq \mu\leq t  \\
         Q_1(z)~g(s)+Q_2(z)~g(t)=u,
    \end{array} \right.
\end{equation}
Let $A(z):G\to L_2^n([s,t])\oplus \mathbb{C}^n$, the corresponding operator with parameter $z\in\mathcal{G}$ will be of the form
\begin{equation}
A(z)g:=
     \left[
    \begin{array}{ll}
        g'-P(z)~g  \\
         Q_1(z)~g(s)+Q_2(z)~g(t)
    \end{array} \right].
\end{equation}
$A(.)$ is a Fredholm operator valued holomorphic function on $\mathcal{G}$ with index zero at every point.
The boundary value problem considered in (\ref{BP}) is called the boundary eigenvalue problem. The details of this kind of problem can be found in \cite{Mennicken}.

\begin{examp}
Keeping all the notations same as in the above discussion we consider $\mathcal{G}=\mathbb{C}$ and 
$$E_+(z)g:=
     \left[
    \begin{array}{ll}
        g'-P_+(z)~g  \\
         Q_1^+(z)~g(s)+Q_2^+(z)~g(t)
    \end{array} \right]$$
    with the fact that, there exists $z_+\in\mathbb{C}$ such that 
    $$Q_1^+(z_+)=I_n\hspace{.5cm}\mbox{and}\hspace{.5cm}Q_2^+(z_+)=0_n,$$
    where $I_n$ is the $n\times n$ identity matrix and $0_n$ is the $n\times n$ zero matrix. Also consider
    $$E_-(z)g:=
     \left[
    \begin{array}{ll}
        g'-P_-(z)~g  \\
         Q_1^-(z)~g(s)+Q_2^-(z)~g(t)
    \end{array} \right]$$
    with the fact that, there exists $z_-\in\mathbb{C}$ such that
    $$Q_1^-(z_-)=I_n\hspace{.5cm}\mbox{and}\hspace{.5cm}Q_2^-(z_-)=0_n.$$
    Now it is easy to observe that $E_+(z_+),E_-(z_-):G\to L_2^n([s,t])\oplus \mathbb{C}^n$ both are invertible and the corresponding inverse operators are
    $$E_+(z_+)^{-1}\begin{pmatrix}
    g\\
    v
    \end{pmatrix}(\mu)=e^{(\mu-s)P_+(z_+)}v+\int_s^\mu e^{(\mu-x)P_+(z_+)}g(x)~dx$$
    and
    $$E_-(z_-)^{-1}\begin{pmatrix}
    g\\
    v
    \end{pmatrix}(\mu)=e^{(\mu-s)P_-(z_-)}v+\int_s^\mu e^{(\mu-x)P_-(z_-)}g(x)~dx$$
    where $g\in L_2^n([s,t])$ and $v\in \mathbb{C}^n$. Now at this point, if we have the $B(G)$- valued function $F:=E_+^{-1}E_-$ belongs to $\mathcal{S}^{in}\cap\mathcal{S}^{in}_*$, then the pair of operator valued functions $$\mathfrak{E}(z)=(E_-(z),E_+(z)) \hspace{.5cm}\mbox{for every}~z\in\mathbb{C}$$
    will represent a de Branges operator.
\end{examp}
We construct the next example with the help of linear operator pencils. More about linear operator pencils can be found in \cite{Gohberg1}, \cite{Mennicken}. Let $A,B\in B(\mathfrak{X})$, then the linear operator pencil
$$S(z)=A-zB$$
is a $B(\mathfrak{X})$-valued entire function. Now suppose for some $z_0\in\mathbb{C}$, $S(z_0)$ is invertible, then we can express $S(z)$ as
\begin{equation}
\label{LP}
S(z)=A-zB=(z_0-z)(A-z_0B)\left[(z_0-z)^{-1}I+(A-z_0B)^{-1}B\right].
\end{equation}
Now along with the invertible condition, if we choose $A\in B(\mathfrak{X})$ and $B$ is compact, it is clear from (\ref{LP}) that $S(z)$ is a Fredholm operator valued entire function with index zero.
\begin{examp}
Suppose $E_+(z)=A-zB$ and $E_-(z)=C-zD$, where $A, C\in B(\mathfrak{X})$ and $B, D$ are compact operators. Also, $E_+$ and $E_-$ both are invertible at least at one point in $\mathbb{C}$. If $E_-^*(z)E_-(z)=I$ for all $z\in\mathbb{C}$ and $E_+(.)^{-1}\in \mathcal{S}^{in}\cap\mathcal{S}^{in}_*$, then the pair of operator valued functions 
$$\mathfrak{E}(z)=(E_-(z),E_+(z)) \hspace{.5cm}\mbox{for every}~z\in\mathbb{C}$$
will represent a de Branges operator.
\end{examp}

The following example involves the system of differential equations of operator valued functions.
\begin{examp}
Let us consider the following initial value problem
\begin{equation}
\label{deq}
\frac{d F_r(z)}{d r}=izF_r(z)j_H+F_r(z)Q(r)\hspace{.4cm}0\leq r\leq a,z\in\mathbb{C}
\end{equation}
with the initial condition (given in the matrix form)
\begin{equation}
F_0(z)=[I_H\hspace{.3cm}I_H],
\end{equation}
where $H$ is a complex seperable Hilbert space, 
$$F_r(z)=[E_-^r(z)\hspace{.3cm}E_+^r(z)]:H\oplus H\to H$$
and
$$j_H=\left[
    \begin{array}{ll}
        I_H\hspace{.5cm}0  \\
         0\hspace{.3cm}-I_H
    \end{array} \right].$$
    Also 
    $$Q(r)=\left[
    \begin{array}{ll}
        0\hspace{.5cm}q(r)  \\
         q(r)^*\hspace{.3cm}0
    \end{array} \right]:H\oplus H\to H\oplus H,$$
    Where $q(r)$ is a $B(H)$-valued continuous function on $[0,a]$.
    Suppose the solution $F_r(z)$ of (\ref{deq}) has the property that both $E_-^r(z)$ and $ E_+^r(z)$ are Fredholm operator valued entire functions for all $0\leq r\leq a$.
 Now for every $\xi\in\mathbb{C}$ it can be easily shown that
\begin{equation}
\frac{d}{ds}\{F_s(z)j_HF_s(\xi)^*\}=i(z-\overline{\xi})F_s(z)F_s(\xi)^*.\label{ex2}
\end{equation}
Observe that $F_0(z)j_HF_0(\xi)^*=0$. Now integrating both sides of (\ref{ex2}) from $0$ to $r$, we get
\begin{equation}
F_r(z)j_HF_r(\xi)^*=i(z-\overline{\xi})\int_0^r F_s(z)F_s(\xi)^*ds.\label{ex3}
\end{equation}
Now if we use the matrix form of $F_r(z)$ in (\ref{ex3}), we get
\begin{equation}
\frac{E_+^r(z)E_+^r(\xi)^*-E_-^r(z)E_-^r(\xi)^*}{\rho_\xi(z)}=\frac{1}{2\pi}\int_0^r F_s(z)F_s(\xi)^*ds.
\end{equation}
Now if we have $\int_0^r F_s(\xi)F_s(\xi)^*ds\succ 0$ for $\xi\in\mathbb{C}_+\cup \mathbb{C}_-$ and there exists $\xi_0\in\mathbb{C}_+$ such that $\int_0^r F_s(\xi_0)F_s(\xi_0)^*ds$, $\int_0^r F_s(\overline{\xi_0})F_s(\overline{\xi_0})^*ds$ both are invertible and 
 $E_-^r(\overline{\xi_0})$, $E_+^r(\xi_0)$ both are self adjoint, then the pair of operator valued functions $\mathfrak{E}^r(z)=(E_-^r(z),~E_+^r(z))$ will represent de Branges operator for all $r\in [0,a]$. 
\end{examp}
For a clear explanation of the last example, see the converse part of the Theorem \ref{de Branges charac}.
\begin{rmk}
The system of differential equations that appeared in the last example was studied widely in the literature. For example, see \cite{Krein system}, where the functions under consideration were scalar and matrix valued.
\end{rmk}

\section{Condition for $\mathcal{H}_\beta$ and $\mathcal{H}_{\overline{\beta}}$ to be isometrically isomorphic}
\label{Se 7}
In this section, we consider an RKHS $\mathcal{H}$ of $\mathfrak{X}$-valued entire functions with RK $K_\xi(z)$ and for some $\beta\in\mathbb{C}$
$$\mathcal{H}_\beta=\{g\in\mathcal{H}:g(\beta)=0\}.$$
It is clear that $\mathcal{H}_\beta$ is a closed subspace of $\mathcal{H}$ and thus an RKHS. Our main goal is to find a condition such that $\mathcal{H}_\beta$ and $\mathcal{H}_{\overline{\beta}}$ are isometrically isomorphic for some $\beta\in\mathbb{C}_+$. Recall that $\mathfrak{T}$ is the multiplication operator in $\mathcal{H}$ where the multiplication is by the independent variable with domain $\mathcal{D}$, which is also a closed operator. The idea of the Moore-Penrose inverse of a bounded linear operator with closed range will be used in this section, and details about it can be found in \cite{Wang}. Some of the results of this section can be found in \cite{ArD20}, where the operator $K_\beta(\beta)$ is invertible.
\begin{lemma}
\label{Kernel}
Suppose $\mathcal{H}$ is an RKHS of $\mathfrak{X}$-valued entire functions with RK $K_\xi(z)$ and assume that $K_\beta(\beta)$ has closed range  for some $\beta\in\mathbb{C}$. Then the following assertions are true:
\begin{enumerate}
\item The RK of $\mathcal{H}_\beta$ can be expressed as
\begin{equation}
K_\xi^\beta(z)=K_\xi(z)-K_\beta(z)K_\beta(\beta)^\dagger K_\xi(\beta),\label{13}
\end{equation}
where $K_\beta(\beta)^\dagger$ is the Moore-Penrose inverse of the operator $K_\beta(\beta)$. 
\item If $\Pi_\beta$ is the orthogonal projection of $\mathcal{H}$ onto $\mathcal{H}_\beta^\perp$ then 
\begin{equation}
\Pi_\beta (g)=K_\beta K_\beta(\beta)^\dagger g(\beta)\hspace{.5cm}\mbox{for all}~g\in\mathcal{H}\label{14}
\end{equation}
and 
\begin{equation}
\mathcal{H}_\beta^\perp=\{K_\beta u:u\in\mathfrak{X}\}=\{K_\beta K_\beta(\beta)^\dagger v:v\in\mathfrak{X}\}.\label{15}
\end{equation}
\item The following equivalence condition hold:
\begin{equation}
R_\beta \mathcal{H}_\beta\subseteq \mathcal{H}~ \mbox{if and only if}~ R_\beta\mathcal{H}_\beta=\mathcal{D}.\label{16}
\end{equation}
\item $\mbox{rng}K_\beta(\beta)=\mbox{rng}(\delta_\beta)$.
\end{enumerate}
Moreover, if $K_\beta(\beta)$ is invertible, then $K_\beta(\beta)^\dagger$ can be replaced by $K_\beta(\beta)^{-1}$ in (\ref{13}) and (\ref{14}).
\end{lemma}
\begin{proof}
For details of the proof see Lemma $2.6$ in \cite{JFA}.
\end{proof}
\begin{lemma}
\label{Unitary map}
If in the setting of lemma $7.1$, the equivalence condition in (\ref{16}) also holds, then the following implications are true:
\begin{enumerate}
\item $R_\beta\in B(\mathcal{H}_\beta,\mathcal{H})$.
\item $\mathcal{H}_\beta=\mbox{rng}(\mathfrak{T}-\beta I)$ and 
\begin{equation}
\mbox{rng}(\mathfrak{T}-\beta I)^\perp =\{K_\beta u:u\in\mathfrak{X}\}=\{K_\beta K_\beta(\beta)^\dagger u:u\in\mathfrak{X}\}.\label{19}
\end{equation}
\item If $K_\beta(\beta)\succ 0$, then there exists a unitary operator $T_\beta$ from $\mbox{rng}(\mathfrak{T}-\beta I)^\perp$ to $\mbox{rng}K_\beta(\beta)$.
\end{enumerate} 
\end{lemma}
\begin{proof}
It is clear that $R_\beta$ is linear. To verify $(1)$, we first show that $R_\beta$ is a closed operator, and the rest of the arguments will be clear from the closed graph theorem. Let $\{g_n:n\in\mathbb{N}\}\subseteq \mathcal{H}_\beta$ be such that $g_n\to g$ and $R_\beta g_n\to h$ as $n\to\infty$. Then $g\in\mathcal{H}_\beta$ as in RKHS norm convergence implies pointwise convergence. Also, for $\xi\neq \beta$,
$$\frac{g(\xi)-g(\beta)}{\xi-\beta}=\lim_{n\to\infty}\frac{g_n(\xi)-g_n(\beta)}{\xi-\beta}=\lim_{n\to\infty}R_\beta g_n(\xi)=h(\xi).$$
Thus $h(\xi)=R_\beta g(\xi)$ for all $\xi\in\mathbb{C}$ as $h$ and $R_\beta g$ are entire functions. This implies that the operator $R_\beta$ is closed.\\
For every $g\in\mathcal{H}_\beta$,
\begin{equation}
(\mathfrak{T}-\beta I)R_\beta g=g,\label{21}
\end{equation}
as for $\xi\neq \beta$ the following holds
$$(\mathfrak{T}-\beta I)R_\beta g(\xi)=\xi\left[\frac{g(\xi)-g(\beta)}{\xi-\beta}\right]-\beta\left[\frac{g(\xi)-g(\beta)}{\xi-\beta}\right]=g(\xi).$$
 This implies
 \begin{equation}
 \mbox{rng}(\mathfrak{T}-\beta I)=\{(\mathfrak{T}-\beta I)g:g\in\mathcal{D}\}=\{(\mathfrak{T}-\beta I)R_\beta f:f\in\mathcal{H}_\beta\}=\mathcal{H}_\beta.
 \end{equation}
 Thus the first assertion of $(2)$ holds. Also from the preceding lemma (\ref{19}) is straightforward. \\
 Since $K_\beta(\beta)\succ 0$ and has closed range, $\mbox{rng}K_\beta(\beta)=\mbox{rng}K_\beta(\beta)^{\frac{1}{2}}$ (see \cite{Fillmore}). Now, in view of (\ref{19}), we consider a map $T_\beta:\mbox{rng}(\mathfrak{T}-\beta I)^\perp\to \mbox{rng}K_\beta(\beta)$ defined by
 \begin{equation}
 T_\beta(K_\beta u)=K_\beta(\beta)^{\frac{1}{2}}u\hspace{.5cm}\mbox{for all}~ u\in \mathfrak{X}.
 \end{equation}
 It is clear that $T_\beta$ is linear and bijective and for any $u\in\mathfrak{X}$
 \begin{align*}
 ||K_\beta u||_\mathcal{H}^2=\langle K_\beta u,K_\beta u\rangle_\mathcal{H}& =\langle K_\beta(\beta)u,u\rangle_\mathfrak{X}\\
 & =\langle K_\beta(\beta)^{\frac{1}{2}}u,K_\beta(\beta)^{\frac{1}{2}}u\rangle_\mathfrak{X}\\
 & =||K_\beta(\beta)^{\frac{1}{2}}u||_\mathfrak{X}^2.
 \end{align*}
 Hence $T_\beta$ is a unitary operator.
\end{proof}

\begin{lemma}
In addition to the setting of Lemma $7.2$, if we assume $\mathcal{D}$ to be dense in $\mathcal{H}$, then for all $u\in\mathfrak{X}$, $K_\beta u$ are the eigen functions of $\mathfrak{T}^*$ corresponding to the eigenvalue $\overline{\beta}$. Also, if we assume $\mathfrak{T}$ to be symmetric, then for some $u\in\mathfrak{X}$, the following equivalence condition holds
\begin{equation}
\rho_\beta K_\beta u\in\mathcal{H}~\mbox{iff}~ K_\beta u=0.\label{20}
\end{equation}
\end{lemma}
\begin{proof}
Let $g\in\mathcal{D}$ and $u\in\mathfrak{X}$, then
	$$\langle \mathfrak{T} g,K_\beta u\rangle_\mathcal{H} = \langle (\mathfrak{T} g)(\beta),u\rangle_\mathfrak{X} = \langle \beta g(\beta),u\rangle_\mathfrak{X} = \langle g(\beta),\overline{\beta}u\rangle_\mathfrak{X} = \langle g,\overline{\beta}K_\beta u\rangle_\mathfrak{X}.$$
	Therefore, $\langle \mathfrak{T} g,K_\beta u\rangle_\mathcal{H}=\langle g,\overline{\beta}K_\beta u\rangle_\mathfrak{X}$ for all $g\in\mathcal{D}$, which proves the first assertion.\\
	To show the equivalence condition in (\ref{20}), it is sufficient to show that for some $u\in\mathfrak{X}$, $\rho_\beta K_\beta u\in\mathcal{H}$ implies $K_\beta u=0$, as the opposite direction is self-evident. Now $\rho_\beta K_\beta u\in\mathcal{H}_{\overline{\beta}}$ and $R_{\overline{\beta}}\rho_\beta K_\beta u=-2\pi iK_\beta u\in\mathcal{H}$, which implies $K_\beta u\in\mathcal{D}$. Let $g\in\mathcal{D}$ and $g=R_\beta h$ for some $h\in\mathcal{H}_\beta$ and observe that
	\begin{equation}
	\rho_\beta K_\beta u=-2\pi i(\mathfrak{T}-\overline{\beta}I)K_\beta u.
	\end{equation}
	Then,
	\begin{align*}
	\langle g,\rho_\beta K_\beta u\rangle_\mathcal{H} & =2\pi i\langle g,(\mathfrak{T}-\overline{\beta}I)K_\beta u\rangle_\mathcal{H}\\
	& =2\pi i\langle (\mathfrak{T}-\beta I)R_\beta h,K_\beta u\rangle_\mathcal{H}\\
	& =2\pi i\langle h,K_\beta u\rangle_\mathcal{H}\\
	& =2\pi i\langle h(\beta),u\rangle_\mathfrak{X}\\
	& =0.
	\end{align*}
	Therefore, $\rho_\beta K_\beta u\perp\mathcal{D}$. Now, due to the additional density condition on $\mathcal{D}$, it is clear that $\rho_\beta K_\beta u=0$. Thus $\rho_\beta K_\beta u\in\mathcal{H} \Rightarrow K_\beta u=0$.
\end{proof}
The next lemma characterizes one of the conditions mentioned by de Branges in \cite{Branges 4} in terms of the symmetric condition of $\mathfrak{T}$. 
\begin{lemma}
\label{Isomorphism}
Let $\mathcal{H}$ be a nonzero RKHS of $\mathfrak{X}$-valued entire functions with RK $K_\xi(z)$. If for some $\beta\in\mathbb{C}_+$, $K_\beta(\beta),K_{\overline{\beta}}(\overline{\beta})$ have closed range and $R_\beta \mathcal{H}_\beta\subseteq \mathcal{H},R_{\overline{\beta}}\mathcal{H}_{\overline{\beta}}\subseteq \mathcal{H}$. Then
\begin{equation}
(\mathfrak{T}-\overline{\beta}I)R_\beta:\mathcal{H}_\beta\to \mathcal{H}_{\overline{\beta}}
\end{equation}
is an isometrically isomorphism iff the operator $\mathfrak{T}$ is symmetric on $\mathcal{D}$.\\
In particular, if $K_\beta(\beta)$ and $K_{\overline{\beta}}(\overline{\beta})$ are invertible, then also the above equivalence holds.
\end{lemma}
\begin{proof}
Observe that 
\begin{equation}
(\mathfrak{T}-\overline{\beta}I)R_\beta(\mathfrak{T}-\beta I)R_{\overline{\beta}}=I_{\mathcal{H}_{\overline{\beta}}}~~\mbox{and}~~(\mathfrak{T}-\beta I)R_{\overline{\beta}}(\mathfrak{T}-\overline{\beta}I)R_\beta=I_{\mathcal{H}_\beta}.
\end{equation}
Now to prove the lemma, we only need to show that $(\mathfrak{T}-\overline{\beta}I)R_\beta$ satisfies the norm preserving property if and only if $\mathfrak{T}$ is symmetric, i.e., for all $g\in\mathcal{H}_\beta$
\begin{equation}
||(\mathfrak{T}-\overline{\beta}I)R_\beta g||_\mathcal{H}=||g||_\mathcal{H}\label{22}
\end{equation}
if and only if $\mathfrak{T}$ is symmetric. In view of (\ref{21}) and by using polarization identity, it is easy to observe that (\ref{22}) holds if and only if
\begin{equation}
\langle (\mathfrak{T}-\overline{\beta}I)R_\beta g,(\mathfrak{T}-\overline{\beta}I)R_\beta h\rangle_\mathcal{H}=\langle g,h\rangle_\mathcal{H}\label{23}
\end{equation}
for every $g,h\in\mathcal{H}_\beta$. Also, for every $g,h\in\mathcal{H}_\beta$, the following evaluation
\begin{align*}
 \langle (\mathfrak{T}-\overline{\beta} I)R_\beta g,(\mathfrak{T}-\overline{\beta} I)R_\beta h\rangle_\mathcal{H}& =\langle g,h\rangle_\mathcal{H} \\
\Updownarrow\\
 \langle \{(\mathfrak{T}-\beta I)+(\beta-\overline{\beta})I\}R_\beta g,\{&(\mathfrak{T}-\beta I)+(\beta-\overline{\beta})I\}R_\beta  h\rangle_\mathcal{H} =\langle g,h\rangle_\mathcal{H} \\
	\Updownarrow\\
 \langle (\mathfrak{T}-{\beta} I)R_\beta g,(\mathfrak{T}-{\beta} I)R_\beta h\rangle_\mathcal{H}&+\langle (\mathfrak{T}-{\beta} I)R_\beta g,(\beta-\overline{\beta})R_\beta h\rangle_\mathcal{H}\\
 +\langle (\beta-\overline{\beta})R_\beta g,(\mathfrak{T}-{\beta} I)R_\beta h\rangle_\mathcal{H}
	& +\langle (\beta-\overline{\beta})R_\beta g,(\beta-\overline{\beta})R_\beta h\rangle_\mathcal{H} =\langle g,h\rangle_\mathcal{H} 	\\
	\Updownarrow\\
 (\overline{\beta}-\beta)[\langle (\mathfrak{T}-{\beta} I)R_\beta g,R_\beta h\rangle_\mathcal{H}&-\langle R_\beta g,(\mathfrak{T}-{\beta} I)R_\beta h\rangle_\mathcal{H}\\
 &+(\beta-\overline{\beta})\langle R_\beta g,R_\beta h\rangle_\mathcal{H}] =0\\
	\end{align*}
	guarantees that (\ref{23}) holds if and only if
	\begin{equation}
	\langle (\mathfrak{T}-{\beta} I)R_\beta g,R_\beta h\rangle_\mathcal{H}-\langle R_\beta g,(\mathfrak{T}-{\beta} I)R_\beta h\rangle_\mathcal{H}+(\beta-\overline{\beta})\langle R_\beta g,R_\beta h\rangle_\mathcal{H}=0.\label{24}
	\end{equation}
	Again for every $g,h\in\mathcal{H}_\beta$, (\ref{24}) holds if and only if
	\begin{equation}
	\langle \mathfrak{T}R_\beta g,R_\beta h\rangle_\mathcal{H}=\langle R_\beta g,\mathfrak{T}R_\beta h\rangle_\mathcal{H}.\label{25}
	\end{equation}
	The following evaluation
	\begin{align*}
	\langle (\mathfrak{T}-{\beta} I)R_\beta g,R_\beta h\rangle_\mathcal{H}-&\langle R_\beta g,(\mathfrak{T}-{\beta} I)R_\beta h\rangle_\mathcal{H}+(\beta-\overline{\beta})\langle R_\beta g,R_\beta h\rangle_\mathcal{H}=0\\
	\Updownarrow\\
	 \langle \{(\mathfrak{T}-{\beta} I)+(\beta-\overline{\beta})I\}&R_\beta g,R_\beta h\rangle_\mathcal{H}-\langle R_\beta g,(\mathfrak{T}-{\beta} I)R_\beta h\rangle_\mathcal{H}=0\\	
	\Updownarrow\\
	 \langle (\mathfrak{T}-\overline{\beta}I)R_\beta g,R_\beta h\rangle_\mathcal{H}-& \langle R_\beta g,(\mathfrak{T}-{\beta} I)R_\beta h\rangle_\mathcal{H}=0\\
	\Updownarrow\\
	 \langle \mathfrak{T}\hspace{.1cm}R_\beta g,R_\beta h\rangle_\mathcal{H}-\overline{\beta}\langle R_\beta& g,R_\beta h\rangle_\mathcal{H}-\langle R_\beta g,\mathfrak{T}\hspace{.1cm}R_\beta h\rangle_\mathcal{H}+\overline{\beta}\langle R_\beta g,R_\beta h\rangle_\mathcal{H}=0\\
	\Updownarrow\\
	 \langle \mathfrak{T}\hspace{.1cm}R_\beta g,R_\beta h\rangle_\mathcal{H}-\langle R_\beta g,&\mathfrak{T}\hspace{.1cm}R_\beta h\rangle_\mathcal{H}=0\\
	\end{align*}
	proves the above equivalence condition. Since $R_\beta\mathcal{H}_\beta=\mathcal{D}$, the first part of the lemma is proved. The case when $K_\beta(\beta)$ and $K_{\overline{\beta}}(\overline{\beta})$ are invertible can be proved similarly.
\end{proof}

\section{A characterization of the RKHS $\mathcal{B}_\beta(\mathfrak{E})$ }
\label{Se 8}
In this section, we discuss about a characterization of $\mathcal{B}_\beta(\mathfrak{E})$ which was initially given by de Branges for RKHS of scalar valued entire functions (see \cite{Branges 4}). This characterization for the RKHS with $p\times p$ entire matrix valued RK can be found in \cite{JFA}. Our observation is in a more general setting where the RK's are operator valued functions.
\begin{lemma}
\label{Invariance condition}
Let $\mathcal{H}=\mathcal{B}(\mathfrak{E})$ be an RKHS based on a de Branges  operator $\mathfrak{E}(z)=(E_-(z),E_+(z))$ as mentioned in Section~\ref{Se 5}. Then $R_\beta\mathcal{H}_\beta\subseteq \mathcal{H}$ if
\begin{enumerate}
\item $\beta\in\overline{\mathbb{C}_+}$ and $E_+(\beta)$ is an invertible operator.
\item $\beta\in\overline{\mathbb{C}_-}$ and $E_-(\beta)$ is an invertible operator.
\end{enumerate}
\end{lemma}
\begin{proof}
The proof is similar to Lemma $6.4$ in \cite{JFA}.
\end{proof}
\begin{thm}
\label{de Branges charac}
Let $\mathcal{H}$ be an RKHS of $\mathfrak{X}$-valued entire functions with $B(\mathfrak{X})$-valued RK $K_\xi(z)$ and suppose $\beta\in\mathbb{C}_+$ be such that
\begin{equation}
K_\beta(z),K_{\overline{\beta}}(z)\in\Phi(\mathfrak{X})\hspace{.5cm}\mbox{for all}~z\in\mathbb{C}\label{26 a}
\end{equation}
and
\begin{equation}
K_\beta(\beta),K_{\overline{\beta}}(\overline{\beta})~\mbox{are invertible}.\label{26}
\end{equation}
Then the RKHS $\mathcal{H}$ is same as the de Branges space $\mathcal{B}_\beta(\mathfrak{E})$ iff 
\begin{equation}
R_\beta\mathcal{H}_\beta\subseteq \mathcal{H},\hspace{.5cm}R_{\overline{\beta}}\mathcal{H}_{\overline{\beta}}\subseteq \mathcal{H}\label{27}
\end{equation}
and 
\begin{equation}
(\mathfrak{T}-\overline{\beta}I)R_\beta:\mathcal{H}_\beta\to \mathcal{H}_{\overline{\beta}}\label{28}
\end{equation}
is an isometrically isomorphism.
\end{thm}
\begin{proof}
The proof will be similar with few exceptions to Theorem $7.1$ in \cite{JFA}. So here we mostly avoid similar evaluations. Since $K_\beta(\beta)$ and $K_{\overline{\beta}}(\overline{\beta})$ are invertible, in view of Lemma \ref{Kernel}, the reproducing kernels of $\mathcal{H}_\beta$ and $\mathcal{H}_{\overline{\beta}}$ are
\begin{equation}
K_\xi^\beta(z)=K_\xi(z)-K_\beta(z)K_\beta(\beta)^{-1} K_\xi(\beta)\label{29}
\end{equation}
and
\begin{equation}
K_\xi^{\overline{\beta}}(z)=K_\xi(z)-K_{\overline{\beta}}(z)K_{\overline{\beta}}(\overline{\beta})^{-1} K_\xi(\overline{\beta})
\end{equation}
respectively. Also, for any $g\in\mathcal{H}_\beta$ and $z\neq \beta$
\begin{equation}
((\mathfrak{T}-\overline{\beta}I)R_\beta g)(z)=\frac{z-\overline{\beta}}{z-\beta}g(z).\label{30}
\end{equation}
First, suppose $\mathcal{H}$ satisfies the constraints in (\ref{27}) and (\ref{28}). Then
\begin{equation}
\frac{z-\overline{\beta}}{z-\beta}K_\xi^\beta(z)=\frac{\overline{\xi}-\overline{\beta}}{\overline{\xi}-\beta}K_\xi^{\overline{\beta}}(z).
\end{equation}
Now consider 
\begin{equation}
E_+(z)=\rho_\beta(z)\rho_\beta(\beta)^{-\frac{1}{2}}K_\beta(z)K_\beta(\beta)^{-\frac{1}{2}}\label{33}
\end{equation}
and
\begin{equation}
E_-(z)=-\rho_{\overline{\beta}}(z)\rho_\beta(\beta)^{-\frac{1}{2}}K_{\overline{\beta}}(z)K_{\overline{\beta}}(\overline{\beta})^{-\frac{1}{2}}.\label{34}
\end{equation}
Then in view of (\ref{26 a}), $E_+(z)$, $E_-(z)$ are entire and belong to $\Phi(\mathfrak{X})$ for all $z\in\mathbb{C}$. Also  $E_+(\beta)=\rho_\beta(\beta)^\frac{1}{2}K_\beta(\beta)^\frac{1}{2}$ and $E_-(\overline{\beta})=\rho_\beta(\beta)^\frac{1}{2}K_{\overline{\beta}}(\overline{\beta})^\frac{1}{2}$. Thus $E_+(\beta)$ and $E_-(\overline{\beta})$ both are invertible and selfadjoint. Also for $z\neq \overline{\xi}$
\begin{equation}
K_\xi(z)= \frac{E_+(z)E_+(\xi)^*-E_-(z)E_-(\xi)^*}{\rho_\xi(z)} . 
\end{equation}
Therefore,
$$E_+(\xi)E_+(\xi)^*-E_-(\xi)E_-(\xi)^*=\rho_\xi(\xi)K_\xi(\xi)\succeq 0$$
for $\xi\in\mathbb{C}_+$ and 
\begin{equation}
E_+(\xi)E_+(\xi)^*-E_-(\xi)E_-(\xi)^*=0\label{35}
\end{equation}
 for $\xi\in\mathbb{R}$. Thus, $E_+^{-1}E_-\in \mathcal{S}^{in}\cap\mathcal{S}^{in}_*$ and the corresponding pair of operator valued functions $\mathfrak{E}(z)=(E_-(z),E_+(z))$ is a de Branges operator. Since the RK's of the spaces $\mathcal{H}$ and  $\mathcal{B}_\beta(\mathfrak{E})$ are equal, $\mathcal{H}=\mathcal{B}_\beta(\mathfrak{E})$.\\
Conversely, let us assume that $\mathcal{H}=\mathcal{B}_\beta(\mathfrak{E})$ and (\ref{26 a}), (\ref{26}) hold. The constraint in (\ref{26}) gives
\begin{equation}
E_+(\beta)E_+(\beta)^*\succ E_-(\beta)E_-(\beta)^*\hspace{.2cm}\mbox{and}\hspace{.2cm}E_-(\overline{\beta})E_-(\overline{\beta})^*\succ E_+(\overline{\beta})E_+(\overline{\beta})^*.\label{31}
\end{equation}
This implies that $E_+(\beta)^*$ and $E_-(\beta)^*$ both are injective. Also in view of Theorem $1$ in \cite{Douglas}, we have
$$\mbox{rng}E_-(\beta)\subseteq \mbox{rng}E_+(\beta)\hspace{.2cm}\mbox{and}\hspace{.2cm}\mbox{rng}E_+(\overline{\beta})\subseteq\mbox{rng}E_-(\overline{\beta}).$$
Since $K_\beta(\beta)$ and $K_{\overline{\beta}}(\overline{\beta})$ both are invertible $E_+(\beta)$ and $E_-(\overline{\beta})$ both are surjective. Thus $E_+(\beta)$ and $E_-(\overline{\beta})$ both are invertible. Now from the preceding lemma, we have 
$$R_\beta\mathcal{H}_\beta\subseteq \mathcal{H}\hspace{.2cm}\mbox{and}\hspace{.2cm}R_{\overline{\beta}}\mathcal{H}_{\overline{\beta}}\subseteq \mathcal{H}.$$
At this point if we prove the norm preserving condition for the operator $(\mathfrak{T}-\overline{\beta}I)R_\beta$ then the rest of the proof follows from Lemma \ref{Isomorphism}. Suppose $g\in\mathcal{H}_\beta$, then $(\mathfrak{T}-\overline{\beta}I)R_\beta g\in\mathcal{H}_{\overline{\beta}}$ and
\begin{align*}
||(\mathfrak{T}-\overline{\beta}I)R_\beta g||^2_{\mathcal{B}_\beta(\mathfrak{E})}& =\int_{-\infty}^\infty ||\frac{x-\overline{\beta}}{x-\beta}(E_+^{-1}g)(x)||^2dx\\
& =\int_{-\infty}^\infty ||(E_+^{-1}g)(x)||^2dx\\
& =||g||^2_{\mathcal{B}_\beta(\mathfrak{E})}.
\end{align*}
\end{proof}

\section{Connection between $\mathfrak{T}$ and the de Branges spaces}
\label{Se 9}
In the present section, we describe the parametrization and canonical description of selfadjoint extensions of $\mathfrak{T}$ with an arbitrary domain $\mathcal{D}$, using the unitary operator $V:\text{rng}\hspace{.1cm}K_\beta(\beta)\to \text{rng}\hspace{.1cm}K_{\overline{\beta}}(\overline{\beta})$ as a parameter. Then with the help of these selfadjoint extensions, we will see that the de Branges space $\mathcal{B}_\beta(\mathfrak{E})$ has Kramer sampling property. Details about the selfadjoint extension of the multiplication operator can be found in \cite{Akhiezer}, \cite{Yuri book}. In the setting of RKHS, which consists of entire $p\times 1$ vector valued functions, the parametrization and canonical description of selfadjoint extensions of the operator $\mathfrak{T}$ with non dense domain $\mathcal{D}$ can be found in this paper \cite{JFA}, where the parameters are $p\times p$ unitary matrices. 
\begin{thm}
\label{Selfadjoint extension}
Suppose $\mathcal{H}$ is an RKHS of $\mathfrak{X}$-valued entire functions with RK $K_\xi(z)$ having at least one nonzero vector, the operator $\mathfrak{T}$ is assumed to be symmetric in its domain $\mathcal{D}$ and for some $\beta\in\mathbb{C}_+$
\begin{itemize}
\item $K_\beta(\beta)$ has closed range, $K_\beta(\beta)\succ 0$ and $R_\beta\mathcal{H}_\beta\subseteq \mathcal{H}$
\item $K_{\overline{\beta}}(\overline{\beta})$ has closed range, $K_{\overline{\beta}}(\overline{\beta})\succ 0$ and $R_{\overline{\beta}}\mathcal{H}_{\overline{\beta}}\subseteq \mathcal{H}$.
\end{itemize}
Then there exists a unitary operator $V:\mbox{rng}K_\beta(\beta)\to\mbox{rng}K_{\overline{\beta}}(\overline{\beta})$ such that the following implications are true:
\begin{enumerate}
\item The following sum
\begin{equation}
\{(T_\beta^{-1}+T_{\overline{\beta}}^{-1}V)u:u\in\mbox{rng}K_\beta(\beta)\}+\mathcal{D}\label{32}
\end{equation}
is direct, where $T_\beta$ as in Lemma \ref{Unitary map}.
\item The operator $\mathfrak{T}_V$ defined as
\begin{equation}
\mathfrak{T}_V(g+T_\beta^{-1}u+T_{\overline{\beta}}^{-1}Vu)=\mathfrak{T}g+\overline{\beta}T_\beta^{-1}u+\beta T_{\overline{\beta}}^{-1}Vu
\end{equation}
with the domain mentioned in (\ref{32}) is a selfadjoint extension of  $\mathfrak{T}$  and the family
 $$\{\mathfrak{T}_V:V~\mbox{is a unitary operator from}~\mbox{rng}K_\beta(\beta)~\mbox{to}~\mbox{rng}K_{\overline{\beta}}(\overline{\beta})~\mbox{satisfying}~(1)\}$$ is the complete list of selfadjoint extensions of $\mathfrak{T}$.
\end{enumerate}
Moreover, if $\mathcal{D}$ is dense in $\mathcal{H}$, then any unitary operator $V:\mbox{rng}K_\beta(\beta)\to\mbox{rng}K_{\overline{\beta}}(\overline{\beta})$ would satisfy (\ref{32}).
\end{thm}
\begin{proof}
Since this is the obvious generalization of Theorem $5.3$ in \cite{JFA} and the technique of the proof is also similar, we avoid the proof. This paper \cite{Straus} will be helpful for the proof.
\end{proof}
\begin{rmk}
If in the setting of the Theorem \ref{Selfadjoint extension}, we assume that $K_\beta(\beta)$ and $K_{\overline{\beta}}(\overline{\beta})$ are invertible, then the unitary operators mentioned in the theorem belong to $B(\mathfrak{X})$, and the domain of $\mathfrak{T}_V$ will be of the following form
\begin{equation}
\{(T_\beta^{-1}+T_{\overline{\beta}}^{-1}V)u:u\in\mathfrak{X}\}\dotplus\mathcal{D}.\label{38}
\end{equation}
Also, the range of the operator $T_\beta$ will be $\mathfrak{X}$, and the inverse will be of the following form
\begin{equation}
T_\beta^{-1}=K_\beta K_\beta(\beta)^{-\frac{1}{2}}.
\end{equation}
\end{rmk}
\begin{thm}
Suppose $\mathcal{H}$ is an RKHS of $\mathfrak{X}$-valued entire functions with $B(\mathfrak{X})$-valued RK $K_\xi(z)$ having at least one nonzero vector and $\beta\in\mathbb{C}_+$ be such that
\begin{enumerate}
\item $K_\beta(z),K_{\overline{\beta}}(z)\in\Phi(\mathfrak{X})$ for all $z\in\mathbb{C}$ and $K_\beta(\beta)$,$K_{\overline{\beta}}(\overline{\beta})$ are invertible.
\item $R_\beta \mathcal{H}_\beta\subseteq \mathcal{H}$ and $R_{\overline{\beta}}\mathcal{H}_{\overline{\beta}}\subseteq \mathcal{H}$.
\item $\mathfrak{T}:\mathcal{D}\to\mathcal{H}$ is symmetric.
\end{enumerate}
Then $\mathcal{H}=\mathcal{B}_\beta(\mathfrak{E})$, where $E_+(z)$ and $E_-(z)$ are as mentioned in (\ref{33}) and (\ref{34}) respectively. Moreover, if for some $\mu\in\mathbb{R}$
\begin{enumerate}
\item[(4)] $K_\mu(\mu)\succ 0$ and $E_+(\mu)$, $E_-(\mu)$ are selfadjoint.
\end{enumerate}
Then the following implications are true:
\begin{enumerate}
\item[(5)] $R_\mu\mathcal{H}_\mu\subseteq \mathcal{H}$, $K_\mu(\mu)$ is invertible, and the operator 
\begin{equation}
V_\mu=(E_-(\mu))^{-1}E_+(\mu)=E_-(\mu)^*(E_+(\mu)^*)^{-1}~\mbox{is unitary}.
\end{equation}
\item[(6)] $V_\mu$ identifies a selfadjoint extension $\mathfrak{T}_{V_\mu}$ of $\mathfrak{T}$.
\item[(7)] $\{K_\mu u:u\in\mathfrak{X}\}$ is the eigenspace corresponding to the eigenvalue $\mu$ of $\mathfrak{T}_{V_\mu}$.
\end{enumerate}
\end{thm}
\begin{proof}
Under the first three assumptions, $\mathcal{H}=\mathcal{B}_\beta(\mathfrak{E})$ follows from Theorem \ref{de Branges charac}. Now for $\mu\in\mathbb{R}$, we have
\begin{equation}
E_+(\mu)E_+(\mu)^*-E_-(\mu)E_-(\mu)^*=0\label{36}
\end{equation}
and 
\begin{equation}
E_+'(\mu)E_+(\mu)^*-E_-'(\mu)E_-(\mu)^*=-2\pi i~ K_\mu^\mathfrak{E}(\mu).\label{37}
\end{equation}
In view of (\ref{36}) and (\ref{37}) we have $E_+(\mu)^*$, $E_-(\mu)^*$ both are injective. Since $E_+(\mu),E_-(\mu)\in\Phi(\mathfrak{X})$ and selfadjoint, both are invertible. Thus $R_\mu\mathcal{H}_\mu\subseteq\mathcal{H}$ follows from Lemma \ref{Invariance condition} and $V_\mu$ is unitary follows from (\ref{36}). Also $R_\mu\in B(\mathcal{H}_\mu,\mathcal{H})$ and $\mu\in \pi(\mathfrak{T})$. This implies $K_\mu(\mu)$ is invertible.\\
Since $V_\mu$ is a unitary operator on $\mathfrak{X}$, to show that $\mathfrak{T}_{V_\mu}$ is a selfadjoint extension of $\mathfrak{T}$, it is sufficient to show that $V_\mu$ satisfies (\ref{38}). Now for $z\in\mathbb{C}$
\begin{align}
T_\beta^{-1}(z) & = K_\beta(z)K_\beta(\beta)^{-\frac{1}{2}}\nonumber\\
& = \frac{\rho_\beta(z)\rho_\beta(\beta)^{-\frac{1}{2}}K_\beta(z)K_\beta(\beta)^{-\frac{1}{2}}}{\rho_\beta(z)\rho_\beta(\beta)^{-\frac{1}{2}}}\nonumber\\
& = \rho_\beta(\beta)^\frac{1}{2}~\frac{E_+(z)}{\rho_\beta(z)}\label{39}\\ 
& = \rho_\beta(\beta)^\frac{1}{2}\left[\frac{E_+(z)}{\rho_\mu(z)}+\frac{\overline{\beta}-\mu}{z-\mu}\frac{E_+(z)}{\rho_\beta(z)}\right].\label{40}
\end{align}
Similarly,
\begin{align}
T_{\overline{\beta}}^{-1}(z) & = -\rho_\beta(\beta)^\frac{1}{2}~\frac{E_-(z)}{\rho_{\overline{\beta}}(z)}\label{41}\\ 
& = -\rho_\beta(\beta)^\frac{1}{2}\left[\frac{E_-(z)}{\rho_\mu(z)}+\frac{\beta-\mu}{z-\mu}\frac{E_-(z)}{\rho_{\overline{\beta}}(z)}\right].\label{42}
\end{align}
For any $V\in B(\mathfrak{X})$, we consider the following notation
\begin{equation}
\chi_\xi^V(z)=(\overline{\beta}-\xi)T_\beta^{-1}(z)+(\beta-\xi)T_{\overline{\beta}}^{-1}(z)V.\label{43}
\end{equation}
In particular,
\begin{equation}
\chi_\mu^{V_\mu}=\rho_\beta(\beta)^\frac{1}{2}\left[(\overline{\beta}-\mu)\frac{E_+}{\rho_\beta}-(\beta-\mu)\frac{E_-}{\rho_{\overline{\beta}}}\right]
\end{equation}
and $\chi_\mu^{V_\mu}(\mu)=0$. Thus $\mathcal{D}=\{R_\mu\chi_\mu^{V_\mu}u:u\in\mathfrak{X}\}$.\\
Now from the above considerations we have
\begin{equation}
T_\beta^{-1}+T_{\overline{\beta}}^{-1}V_\mu=\rho_\beta(\beta)^\frac{1}{2}\left[\frac{E_+-E_-V_\mu}{\rho_\mu}\right]+R_\mu\chi_\mu^{V_\mu}.\label{44}
\end{equation}
Now multiplying (\ref{44}) by $E_+(\mu)^*$ from the right, we get
\begin{equation}
(T_\beta^{-1}(z)+T_{\overline{\beta}}^{-1}(z)V_\mu)E_+(\mu)^*=\rho_\beta(\beta)^\frac{1}{2}K_\mu^\mathfrak{E}(z)+(R_\mu\chi_\mu^{V_\mu})(z)E_+(\mu)^*.\label{45}
\end{equation}
Now if for some $u\in\mathfrak{X}$, $(T_\beta^{-1}+T_{\overline{\beta}}^{-1}V_\mu)E_+(\mu)^*u\in\mathcal{D}$, the above identity gives $\rho_\mu K_\mu^\mathfrak{E}u\in\mathcal{H}$. Since $K_\mu(\mu)\succ 0$, this implies $u=0$. Thus $(2)$ holds and $\mathfrak{T}_{V_\mu}$ is a selfadjoint extension of $\mathfrak{T}$.\\
From (\ref{45}) it is clear that $K_\mu^\mathfrak{E}u$ belongs to the domain of $\mathfrak{T}_{V_\mu}$ and $$(T_{V_\mu}-\mu I)K_\mu^\mathfrak{E}u=0$$ for all $u\in\mathfrak{X}$. Thus, $(3)$ holds as $K_\mu(\mu)$ is invertible.
\end{proof}
\begin{thm}
Suppose $\mathcal{H}$ is an RKHS of $\mathfrak{X}$-valued entire functions with RK $K_\xi(z)$ having at least one nonzero vector such that (\ref{26 a}), (\ref{26}) hold, and  $\{K_{\mu_i}u_i\}$ is an orthogonal basis of $\mathcal{H}$ for $\mu_1,\mu_2,\ldots\in\mathbb{R}$ and $u_1,u_2,\ldots\in\mathfrak{X}$. Then 
\begin{enumerate}
\item $\mathfrak{T}:\mathcal{D}\to\mathcal{H}$ is symmetric.
\item $\mathcal{H}=\mathcal{B}_\beta(\mathfrak{E})$.
\end{enumerate}
Moreover, if $V\in B(\mathfrak{X})$ is a unitary operator satisfying (\ref{38}), $K_\mu(\mu)\succ 0$ and $E_+(\mu), E_-(\mu)$ both are selfadjoint, then
\begin{enumerate}
\item[(3)]  $\mu\in\mathbb{R}$ is an eigenvalue of $\mathfrak{T}_V$ if and only if 
\begin{equation}
\{E_+(\mu)-E_-(\mu)V\}u=0
\end{equation}
and the corresponding eigenfunction
\begin{equation}
g=\lambda K_\mu^\mathfrak{E}(E_+(\mu)^*)^{-1}u
\end{equation} 
for some nonzero $\lambda\in\mathbb{C}$ and nonzero $u\in\mathfrak{X}$. Also, the geometric multiplicity of the eigenvalue $\mu$ is countably infinite.
\item[(4)] If $E_+(\mu)-E_-(\mu)V$ is invertible, then $(\mathfrak{T}_V-\mu I)$ is a closed operator, and $\mu\not\in\sigma(\mathfrak{T}_V)$.
\item[(5)] If $E_+(z)-E_-(z)V\in\Phi(\mathfrak{X})$ for all $z\in\mathbb{C}$ and invertible at least at one point, then $\mathfrak{T}_V$ has a discrete set of eigenvalues.
\end{enumerate}
\end{thm}
\begin{proof}
$(1),~ (2)$ will follow from Theorem \ref{de Branges charac} in association with Lemma \ref{Isomorphism}, once we show the norm preserving property of the operator $(\mathfrak{T}-\overline{\beta}I)R_\beta:\mathcal{H}_\beta\to \mathcal{H}_{\overline{\beta}}$. Let $g\in \mathcal{H}_\beta$, then
\begin{align*}
||(\mathfrak{T}-\overline{\beta}I)R_\beta g||_\mathcal{H}^2 & =\sum_{i=1}^\infty |\langle (\mathfrak{T}-\overline{\beta}I)R_\beta g,\frac{K_{\mu_i}u_i}{||K_{\mu_i}u_i||_\mathcal{H}}\rangle|^2 \\ & =\sum_{i=1}^\infty |\langle ((\mathfrak{T}-\overline{\beta}I)R_\beta g)(\mu_i),\frac{u_i}{||K_{\mu_i}u_i||_\mathcal{H}}\rangle|^2\\
& = \sum_{i=1}^\infty |\frac{\mu_i-\overline{\beta}}{\mu_i-\beta}|^2~|\langle g,\frac{K_{\mu_i}u_i}{||K_{\mu_i}u_i||_\mathcal{H}}\rangle|^2\\
& = \sum_{i=1}^\infty |\langle g,\frac{K_{\mu_i}u_i}{||K_{\mu_i}u_i||_\mathcal{H}}\rangle|^2=||g||^2_\mathcal{H}.\\
\end{align*}
Now suppose $V$ is a unitary operator satisfying (\ref{38}), $K_\mu(\mu)\succ 0$ and $E_+(\mu)$, $E_-(\mu)$ both are selfadjoint, then $R_\mu\mathcal{H}_\mu=\mathcal{D}$. Let $\mu\in\mathbb{R}$ is an eigenvalue of $\mathfrak{T}_V$. Then there exists a nonzero vector $g=h+(T_\beta^{-1}+T_{\overline{\beta}}^{-1}V)u$ in domain of $\mathfrak{T}_V$, where $h\in\mathcal{D}$ and $u\in\mathfrak{X}$. Thus for all $z\in\mathbb{C}$
\begin{equation}
((\mathfrak{T}_V-\mu I)g)(z)=(z-\mu)h(z)+\chi_\mu^V(z)u=0,
\end{equation}
which gives $\chi_\mu^V(\mu)u=0$, $h(z)=-(R_\mu\chi_\mu^V)(z)u$ and $u\neq 0$. This implies $R_\mu\chi_\mu^V u\in\mathcal{D}$ and 
$$g(z)=-(R_\mu\chi_\mu^V)(z)u+(T_\beta^{-1}+T_{\overline{\beta}}^{-1}V)u.$$
Then by using the fact that $\chi_\mu^V(\mu)u=0$, $g$ can be expressed in the following form
\begin{equation}
g(z)=\rho_\beta(\beta)^\frac{1}{2}\left[\frac{E_+(z)-E_-(z)V}{\rho_\mu(z)}u\right].
\end{equation}
It can also be proved that
\begin{equation}
\chi_\mu^V(\mu)u=0\iff(E_+(\mu)-E_-(\mu)V)u=0.
\end{equation}
Since $E_+(\mu)$ and $E_-(\mu)$ both are invertible, we have
\begin{equation}
(E_+(\mu)-E_-(\mu)V)u=0\iff Vu=E_-(\mu)^*(E_+(\mu)^*)^{-1}u.
\end{equation}
This gives 
$$g(z)=\rho_\mu(\mu)^\frac{1}{2}K_\mu^\mathfrak{E}(z)(E_+(\mu)^*)^{-1}u.$$
To show the converse part of $(3)$, we first observe that, if for some $u\neq 0$, $\{E_+(\mu)-E_-(\mu)V\}u=0$ and $g=\lambda K_\mu^\mathfrak{E}(E_+(\mu)^*)^{-1}u$ then $Vu=V_\mu u$ and $\chi_\mu^V u=\chi_\mu^{V_\mu}u$. This implies
$$g=\rho_\beta(\beta)^\frac{1}{2}K_\mu^\mathfrak{E}(E_+(\mu)^*)^{-1}u=-R_\mu\chi_\mu^V u+(T_\beta^{-1}+T_{\overline{\beta}}^{-1}V)u$$ 
belongs to the domain of $\mathfrak{T}_V$ and
\begin{multline*}
(\mathfrak{T}_V-\mu I)g=-\mathfrak{T}(R_\mu\chi_\mu^V u)+(\overline{\beta}T_\beta^{-1}+\beta T_{\overline{\beta}}^{-1}V)u+\mu R_\mu\chi_\mu^V u-\mu(T_\beta^{-1}+T_{\overline{\beta}}^{-1}V)u\\
=\{-\chi_\mu^V+(\overline{\beta}-\mu)T_\beta^{-1}+(\beta-\mu)T_{\overline{\beta}}^{-1}V\}u=0.
\end{multline*}
Thus, $(3)$ holds.\\ Now suppose $E_+(\mu)-E_-(\mu)V$ is invertible.
Since $\mathfrak{T}_V$ is selfadjoint, the operator $(\mathfrak{T}_V-\mu I)$ is closed. To verify  $\mu\not\in\sigma(\mathfrak{T}_V)$ we need to show that  $(\mathfrak{T}_V-\mu I)^{-1}$ exists and is bounded. $(3)$ implies that $(\mathfrak{T}_V-\mu I)$ is injective on domain of $\mathfrak{T}_V$. Since  $\mathfrak{T}_V$ is selfadjoint and $\mu\in\mathbb{R}$, $\mbox{rng}(\mathfrak{T}_V-\mu I)$ is dense in $\mathcal{H}$. Also it can be proved that $(\mathfrak{T}_V-\mu I)$ is surjective (for a similar proof see Theorem $8.5$ in \cite{JFA}). Now the rest of the arguments follow from closed graph theorem. $(5)$ follows from $(3)$ and Theorem \ref{Fredholm thm}.
\end{proof}
The following theorem gives that under some special conditions, the de Branges space $\mathcal{B}_\beta(\mathfrak{E})$ has the Kramer sampling property.
\begin{thm}
Suppose $\mathcal{H}$ is an RKHS of $\mathfrak{X}$-valued entire functions with RK $K_\xi(z)$ having at least one nonzero vector such that (\ref{26 a}), (\ref{26}) hold, and the operator $\mathfrak{T}:\mathcal{D}\to\mathcal{H}$ is assumed to be symmetric. Then the RKHS $\mathcal{H}$ is the de Branges space $\mathcal{B}_\beta(\mathfrak{E})$.\\
Moreover, if $K_\beta(z),K_{\overline{\beta}}(z)$ are invertible for all $z\in\mathbb{R}$ and there exists a unitary operator $V\in B(\mathfrak{X})$ satisfying (\ref{38}) such that $E_+(z)-E_-(z)V\in\Phi(\mathfrak{X})$ for all $z\in\mathbb{C}$ and invertible at least at one point then $\mathcal{B}_\beta(\mathfrak{E})$ has Kramer sampling property.
\end{thm}
\begin{proof}
$\mathcal{H}=\mathcal{B}_\beta(\mathfrak{E})$ follows from Theorem \ref{de Branges charac}. Suppose $V\in B(\mathfrak{X})$ is a unitary operator satisfying all the conditions mentioned in the statement. Then we can have a selfadjoint extension $\mathfrak{T}_V$ of $\mathfrak{T}$ and the spectrum $\sigma(\mathfrak{T}_V)\subseteq \mathbb{R}$. Now for some $\mu\in\mathbb{R}$, if $E_+(\mu)-E_-(\mu)V$ is invertible, then from the preceding theorem, it is clear that $\mu\not\in\sigma(\mathfrak{T}_V)$. This gives 
$$\sigma(\mathfrak{T}_V)=\{\mu\in\mathbb{R}:E_+(\mu)-E_-(\mu)V~\mbox{is not invertible}\},$$
which is precisely the collection of all eigenvalues of $\mathfrak{T}_V$. Also, $\sigma(\mathfrak{T}_V)$ is a discrete set. Since $K_\beta(z),K_{\overline{\beta}}(z)$ are invertible for all $z\in\mathbb{R}$, $E_+(z),E_-(z)$ both are invertible there. Thus the eigenfunctions are of the form $g=K_\mu^\mathfrak{E}u$, and the eigenspaces are countably infinite. Since $\mathfrak{T}_V$ is selfadjoint, any two eigenfunctions corresponding to different eigenvalues are orthogonal, and the Gram-Schmidt orthogonalization process can be used to make the eigen functions orthogonal corresponding to the same eigenvalue. Also, since $\mathfrak{T}_V$ is selfadjoint, the spectral theorem implies that the collection of eigenfunctions is total in $\mathcal{H}$. This completes the proof.
\end{proof}
\section{Entire operators with infinite deficiency indices}
\label{Se 10}
This section revives a functional model problem regarding entire operators with infinite deficiency indices. M. G. Krein introduced and primarily studied these entire operators and made connections with the multiplication operator in a Hilbert space of analytic functions on $\mathbb{C}$. In the fundamental paper \cite{Krein1}, he showed that an entire operator with arbitrary finite equal deficiency indices $(p,p)$ could be considered as the multiplication operator in a Hilbert space of $\mathbb{C}^p$-valued entire functions. Later in this paper \cite{JFA}, it was observed that this Hilbert space is a de Branges space with $p\times p$ matrix valued RK. Krein also studied the entire operators with infinite deficiency indices (see \cite{Krein}), and a similar connection with the multiplication operator in a Hilbert space of $\mathfrak{X}$-valued entire functions has been mentioned here \cite{Gorbachuk} (Appendix I).\\
Since in Section~\ref{Se 5}, we have constructed the de Branges spaces of $\mathfrak{X}$-valued entire functions. It is a natural question whether these newly constructed de Branges spaces can be considered as the functional model for entire operators with infinite deficiency indices. Assume that $Y$ is an infinite dimensional closed subspace of
 $\mathfrak{X}$. Let $E$ is a densely defined closed, simple, symmetric operator on $\mathfrak{X}$ with infinite deficiency indices. We denote $\rho_Y(E)$ as the collection of all $Y$-regular points of $E$, which is defined by 
\begin{equation}
\rho_Y(E):=\{\xi\in\mathbb{C}: \mathfrak{M}_\xi=\mbox{rng}(E-\xi I)=\overline{\mathfrak{M}_\xi}\hspace{.2cm}\mbox{and}\hspace{.2cm}\mathfrak{X}=\mathfrak{M}_\xi\dotplus  Y\}.\label{10.1}
\end{equation}
It is known that $\rho_Y(E)$  is an open subset of $\mathbb{C}$ and every $\xi\in \rho_Y(E)$ is also a point of regular type for $E$. Because of (\ref{10.1}), it is clear that for every $\xi\in \rho_Y(E)$, there exists the projection operator $P_Y(\xi)$, i.e., for every $f\in\mathfrak{X}$, there exists a unique $g\in \mathcal{D}(E)$, the domain of $E$, such that 
$$f=(E-\xi I)g+P_Y(\xi)f.$$
 Also, for every fixed $f\in\mathfrak{X}$, we can consider a map from $\rho_Y(E)$ to $Y$ defined by $\xi\mapsto P_Y(\xi)f$. We denote these $Y$-valued functions as $f_Y$ for every $f\in\mathfrak{X}$ and are defined as $f_Y(\xi)=P_Y(\xi)f$, also assume $\mathcal{H}:=\{f_Y:f\in\mathfrak{X}\}$. Let $\xi\in \rho_Y(E)$, then
$$\mbox{rng}P_Y(\xi)=Y\hspace{.5cm}\mbox{and}\hspace{.5cm}\mbox{ker}P_Y(\xi)=\mbox{rng}(E-\xi I).$$
Since both range and kernel of the projection operator $P_Y(\xi)$ are closed subspaces of $\mathfrak{X}$, $P_Y(\xi)$ is bounded for all $\xi\in \rho_Y(E)$. Then,

$$\mbox{rng}P_Y(\xi)^*=\mathfrak{X}\ominus \mbox{rng}(E-\xi I)\hspace{.5cm}\mbox{and}\hspace{.5cm}\mbox{ker}P_Y(\xi)^*=Y^\perp.$$
Also, for every $\xi\in \rho_Y(E)$, we can have the operator $\mathcal{T}_Y(\xi)\in B(\mathfrak{X})$, which is defined by
$$\mathcal{T}_Y(\xi):=(E-\xi I)^{-1}(I-P_Y(\xi)).$$
 Now following Krein's definition for entire operators, $E$ is an entire operator if $\rho_Y(E)=\mathbb{C}$, and the functions $f_Y$ are entire. This implies that $P_Y(\xi)$ and $\mathcal{T}_Y(\xi)$ both are $B(\mathfrak{X})$- valued entire functions. More properties of these two functions can be found in \cite{Gorbachuk}.
 \begin{lemma}
 \label{L10.1}
 For any $\xi\in \rho_Y(E)$, the restriction of the projection operator $P_Y(\xi)$ on $\mathfrak{M}_\xi^\perp$ is invertible, i.e., the operator $P_Y(\xi)\vert_{\mathfrak{M}_\xi^\perp}:\mathfrak{M}_\xi^\perp\to Y$ is invertible.
 \end{lemma}
 \begin{proof}
 Suppose $f, g\in \mathfrak{M}_\xi^\perp= \mathfrak{X}\ominus \mbox{rng}(E-\xi I)$ be such that $$P_Y(\xi)f=P_Y(\xi)g=h~(\mbox{say}).$$ 
 Then there exist $f_1, g_1\in \mathfrak{M}_\xi$, such that $f=f_1+h$ and $g=g_1+h$. Since $f_1-g_1\in \mathfrak{M}_\xi$ and $f-g\in\mathfrak{M}_\xi^\perp$, $P_Y(\xi)\vert_{\mathfrak{M}_\xi^\perp}$ is one-one.\\
 Now for any $f\in Y$, we have the unique sum $f=g+h$, where $g\in \mathfrak{M}_\xi$ and $h\in \mathfrak{M}_\xi^\perp$. This implies $P_Y(\xi)\vert_{\mathfrak{M}_\xi^\perp}$ is onto.
 \end{proof}
  Since $E$ is simple, the map  $\Psi :\mathfrak{X}\to \mathcal{H}$ defined by $f\mapsto f_Y$ is injective. Thus $\mathcal{H}$ is a vector space with respect to the point wise addition and scalar multiplication. Consider the inner product in $\mathcal{H}$ defined as 
$$\langle f_Y,g_Y\rangle_\mathcal{H}:=\langle f,g\rangle_\mathfrak{X}\hspace{.5cm}\mbox{for all}~f,g\in\mathfrak{X}.$$ 
It is clear that $\mathcal{H}$ is a Hilbert space with respect to the above inner product, and $\Psi$ is a unitary operator.\\
 Let $f\in\mathcal{D}(E)$ and $g\in\mathfrak{X}$ be such that $g=Ef$. For every $\xi\in\mathbb{C}$ there exists unique $f_\xi'\in\mathcal{D}(E)$ such that 
$$f=(E-\xi I) f_\xi'+f_Y(\xi).$$
This gives
\begin{align*}
g=Ef & = (E-\xi I)f+\xi f\\
& = (E-\xi I)f+\xi\{(E-\xi I) f_\xi'+f_Y(\xi)\}\\
& = (E-\xi I)(f+\xi f_\xi')+\xi f_Y(\xi).
\end{align*}
Because of (\ref{10.1}), it is easy to observe that $g_Y(\xi)=\xi f_Y(\xi)$ for all $\xi\in\mathbb{C}$. Thus the operator $E$ on $\mathfrak{X}$ is unitarily equivalent to the multiplication operator on $\mathcal{H}$.\\
Now for any $\xi\in\mathbb{C}$ and $f_Y\in\mathcal{H}$, we have
$$||f_Y(\xi)||_Y=||P_Y(\xi)f||_Y\leq ||P_Y(\xi)||~||f||_\mathfrak{X}=||P_Y(\xi)||~||f_Y||_\mathcal{H}.$$
Since for all $\xi\in\mathbb{C}$, the projection operators $P_Y(\xi)$ are bounded, the point evaluation linear maps in $\mathcal{H}$ are bounded. This implies that $\mathcal{H}$ is an RKHS with the reproducing kernel
$$K_\xi(z)=\delta_z\delta_\xi^*\hspace{.5cm}\mbox{for all}~\xi,z\in\mathbb{C}.$$
Now, let us observe the range and the kernel of the operator $\delta_z:\mathcal{H}\to Y$ for any $z\in\mathbb{C}$. Let $f_Y\in\mathcal{H}$ be such that $\Psi(f)=f_Y$ for $f\in\mathfrak{X}$. Then
$$\delta_z(f_Y)=f_Y(z)=P_Y(z)f.$$
Thus,
\begin{equation}
\mbox{rng}\delta_z=\mbox{rng}P_Y(z)=Y~~\mbox{and}~~\mbox{ker}\delta_z=\{f_Y=\Psi(f):f\in\mbox{rng}(E-zI)\}.\label{aa}
\end{equation}
  This implies
  \begin{equation}
  \mbox{ker}\delta_z^*=\{0\}~~\mbox{and}~~\mbox{rng}\delta_z^*=\{f_Y=\Psi(f):f\in\mathfrak{X}\ominus\mbox{rng}(E-zI)\}.\label{bb}
  \end{equation}
  At this point, we also recall a few facts regarding the generalized Cayley transform. Suppose $E'$ is a selfadjoint extension of the entire operator $E$ within $\mathfrak{X}$. Then the generalized Cayley transform is defined by
  \begin{equation}
  I+(\xi-z)(E'-\xi I)^{-1}=(E'-zI)(E'-\xi I)^{-1}\hspace{.3cm}\mbox{for all}~\xi,z\in\mathbb{C}\setminus \sigma(E').
  \end{equation}
  It is known that the generalized Cayley transform
  \begin{equation}
  \label{CT}
  I+(\xi-z)(E'-\xi I)^{-1}:\mathfrak{M}_{\overline{z}}^\perp\to\mathfrak{M}_{\overline{\xi}}^\perp
  \end{equation}
  is bijective for all $\xi,z\in\mathbb{C}\setminus \sigma(E')$. More details about the generalized Cayley transform can be found in \cite{Gorbachuk} (Chapter $1$, Section $2$). The following lemma has also been collected from \cite{Gorbachuk}, which will provide a necessary motivation for our final problem.
\begin{lemma}
\label{L10.2}
Suppose $\mathfrak{R}_\xi=(E-\xi I)^{-1}$ for all $\xi\in\mathbb{C}$ for the entire operator $E$. Then, for any two numbers $\xi,z\in\mathbb{C}$ the operator
\begin{equation}
I+(\xi-z)\mathfrak{R}_\xi:\mathfrak{M}_\xi\to\mathfrak{M}_z
\end{equation}
is bijective.
\end{lemma}
\begin{proof}
Let $f\in\mathfrak{M}_\xi$, then there exists $g\in\mathcal{D}(E)$ such that $f=(E-\xi I)g$. Now 
\begin{align*}
[I+(\xi-z)\mathfrak{R}_\xi]f & = [I+(\xi-z)\mathfrak{R}_\xi](E-\xi I)g\\
& = (E-\xi I)g+(\xi-z)g\\
& = (E-zI)g\in\mathfrak{M}_z.
\end{align*}
Since every $\xi\in\mathbb{C}$ is a point of regular type of $E$, the operator $(E-\xi I)$ is injective, and this implies the operator $I+(\xi-z)\mathfrak{R}_\xi$ is also injective for every $\xi,z\in\mathbb{C}$. The operator $I+(\xi-z)\mathfrak{R}_\xi$ is also surjective as for any $g\in\mathfrak{M}_z$ with $g=(E-z I)g'$ for $g'\in\mathcal{D}(E)$, the element $f=(E-\xi I)g'\in\mathfrak{M}_\xi$ is the pre-image of $g$.
\end{proof}
Recall that for $z\in\mathbb{C}$, $R_z$ is the generalized backward-shift operator. Suppose $f\in\mathfrak{X}$, then for any $\xi,z\in\mathbb{C}$ there exists $f_\xi',f_z'\in\mathcal{D}(E)$ such that
$$f=(E-\xi I)f_\xi'+f_Y(\xi)=(E-zI)f_z'+f_Y(z).$$
Now a simple calculation gives
$$f_z'=(E-\xi I)\frac{f_\xi'-f_z'}{\xi-z}+\frac{f_Y(\xi)-f_Y(z)}{\xi-z}.$$
This implies the invariance of $\mathcal{H}$ under $R_z$ for all $z\in\mathbb{C}$.\\
Since the operator $E$ on $\mathfrak{X}$ is symmetric and unitarily equivalent to the multiplication operator on $\mathcal{H}$, then the multiplication operator is also symmetric on $\mathcal{H}$. Finally, we summarise all the results we discussed in this section in terms of a theorem, which will also serve the purpose of answering a problem of functional model of entire operators with infinite deficiency indices.
\begin{thm}
Suppose $\mathfrak{X}$ is a complex separable Hilbert space and $E$ be an entire operator with infinite deficiency indices, producing the direct sum decomposition of $\mathfrak{X}$ as mentioned in (\ref{10.1}). Also, suppose for at least one $\beta\in\mathbb{C}_+$ the following conditions hold:
\begin{enumerate}
\item the dimensions of $\mathfrak{M}_\beta\cap\mathfrak{M}_\xi^\perp$ are finite for all $\xi\in\mathbb{C}_-$ and the dimensions of $\mathfrak{M}_{\overline{\beta}}\cap\mathfrak{M}_\xi^\perp$ are finite for all $\xi\in\mathbb{C}_+$.
\end{enumerate}
Then $E$ is unitarily equivalent to the densely defined multiplication operator in a de Branges space $\mathcal{B}_\beta(\mathfrak{E})$. The space $\mathcal{B}_\beta(\mathfrak{E})$ is also invariant under the generalized backward-shift operator $R_z$ for all $z\in\mathbb{C}$.
\end{thm}
\begin{proof}
We begin the proof by observing the intersection of some related closed subspaces of $\mathfrak{X}$. Due to Lemma $2.1$ in \cite{Shtraus}, we have
\begin{equation}
\mathfrak{M}_z\cap\mathfrak{M}_\xi^\perp=\{0\}\hspace{.4cm}\forall ~z,~\xi:\mbox{Im}~z\cdot\mbox{Im}~\xi>0.
\end{equation}
Moreover the following direct sum decomposition
$$\mathfrak{X}=\mathfrak{M}_z\dotplus\mathfrak{M}_\xi^\perp\hspace{.4cm}\forall~z,~\xi:\mbox{Im}~z\cdot\mbox{Im}~\xi>0$$
holds. Now, since any $a\in\mathbb{R}$ is a point of regular type of $E$ the operator $\tilde{E_a}$,
$$\mathcal{D}(\tilde{E_a})=\mathcal{D}(E)\dotplus\mathfrak{M}_a^\perp,~\tilde{E_a}(f_E+\phi_a)=Ef_E+a~\phi_a,~f_E\in\mathcal{D}(E),~\phi_a\in\mathfrak{M}_a^\perp$$
is selfadjoint (see \cite{Shtraus}). Thus every $z\in\mathbb{C}\setminus\mathbb{R}$ is a regular point of $\tilde{E_a}$. This implies for any $z\in\mathbb{C}\setminus\mathbb{R}$ and $f\in\mathfrak{X}$, there exist unique $f_E\in\mathcal{D}(E)$ and $\phi_a\in\mathfrak{M}_a^\perp$ such that
$$f=(\tilde{E_a}-zI)(f_E+\phi_a)=(E-zI)f_E+(a-z)\phi_a.$$
This gives the following direct sum decomposition
$$\mathfrak{X}=\mathfrak{M}_z\dotplus\mathfrak{M}_a^\perp\hspace{.4cm}\forall~ z\in\mathbb{C}\setminus\mathbb{R}~\mbox{and}~\forall~a\in\mathbb{R}.$$
 Thus the intersection
\begin{equation}
\mathfrak{M}_z\cap\mathfrak{M}_a^\perp=\{0\}\hspace{.4cm}\forall~ z\in\mathbb{C}\setminus\mathbb{R}~\mbox{and}~\forall~a\in\mathbb{R}.
\end{equation}
Also 
\begin{equation}
\mathfrak{M}_a\cap\mathfrak{M}_z^\perp=\{0\}\hspace{.4cm}\forall~ z\in\mathbb{C}\setminus\mathbb{R}~\mbox{and}~\forall~a\in\mathbb{R}.
\end{equation}
As $f\in\mathfrak{M}_a\cap\mathfrak{M}_z^\perp$ implies $f\perp\mathfrak{M}_a^\perp$ and $f\perp\mathfrak{M}_z$. Thus $f\perp(\mathfrak{M}_z\dotplus\mathfrak{M}_a^\perp)=\mathfrak{X}$, which implies $f=0$. It can also be proved that the direct sum decomposition
$$\mathfrak{X}=\mathfrak{M}_a\dotplus\mathfrak{M}_z^\perp\hspace{.4cm}\forall~ z\in\mathbb{C}\setminus\mathbb{R}~\mbox{and}~\forall~a\in\mathbb{R}$$
 holds. Now suppose $E'$ is a selfadjoint extension of $E$ within $\mathfrak{X}$, then due to (\ref{CT}) and Lemma \ref{L10.2}, we have
$$I+(\beta-\xi)(E'-\beta I)^{-1}:\mathfrak{M}_\beta\cap \mathfrak{M}_{\overline{\xi}}^\perp\to\mathfrak{M}_\xi\cap \mathfrak{M}_{\overline{\beta}}^\perp$$
and
$$I+(\overline{\beta}-\xi)(E'-\overline{\beta}I)^{-1}:\mathfrak{M}_{\overline{\beta}}\cap \mathfrak{M}_{\overline{\xi}}^\perp\to\mathfrak{M}_\xi\cap \mathfrak{M}_\beta^\perp$$
are bijective for all $\xi\in\mathbb{C}\setminus \sigma(E')$. These observations together with condition $(1)$ imply that the subspaces $\mathfrak{M}_\xi\cap\mathfrak{M}_\beta^\perp$, $\mathfrak{M}_\beta\cap\mathfrak{M}_\xi^\perp$,  $\mathfrak{M}_\xi\cap\mathfrak{M}_{\overline{\beta}}^\perp$ and $\mathfrak{M}_{\overline{\beta}}\cap\mathfrak{M}_\xi^\perp$ are finite dimensional for all $\xi\in\mathbb{C}$.\\
Now, since $K_\beta(\xi)=\delta_\xi\delta_\beta^*$, we have
$$\dim(\ker \delta_\xi\delta_\beta^*)=\dim(\ker \delta_\beta^*)+\dim(\ker \delta_\xi\cap \mbox{rng}\delta_\beta^*).$$
Due to (\ref{aa}) and (\ref{bb}), we have for any $\xi\in\mathbb{C}$,
$$\dim(\ker K_\beta(\xi))=\dim(\mathfrak{M}_\xi\cap\mathfrak{M}_\beta^\perp)~\mbox{and}~\dim(\ker K_\beta(\xi)^*)=\dim(\mathfrak{M}_\beta\cap\mathfrak{M}_\xi^\perp).$$
 Thus the above observation implies that $K_\beta(\xi)\in \Phi(Y)$ for all $\xi\in\mathbb{C}$. Similarly, it can also be observed that $K_{\overline{\beta}}(\xi)\in \Phi(Y)$ for all $\xi\in\mathbb{C}$. 
Also, Lemma \ref{L10.1} implies that $K_\beta(\beta)$ and $K_{\overline{\beta}}(\overline{\beta})$ both are invertible. The rest of the proof follows from the previous discussions in this section and in association with Theorem \ref{de Branges charac} and Lemma \ref{Isomorphism}.
\end{proof}

\section{Connection with the characteristic function of a contraction operator}
\label{Se 11}
In this section, we construct RKHS of $\mathfrak{X}$-valued analytic  functions using the characteristic function of a completely nonunitary (cnu) contraction operator. The underlying idea is to consider those cnu contraction operators whose characteristic functions are inner and invertible on $\mathbb{D}$. These inner functions are then considered on the upper half plane with the help of the conformal map $C$ and construct RKHS using the same technique mentioned in section~\ref{Se 4}. Here we dealt with two situations which will be discussed separately.  Most of the standard results and notations used in this section can be found in \cite{Nagy}.\\
 Let $A\in B(\mathfrak{X})$ be a completely nonunitary contraction operator. Recall that the characteristic function of $A$ is given by
$$C_A(z)=\left.\left[-A+z(I-AA^*)^\frac{1}{2}(I-zA^*)^{-1}(I-A^*A)^\frac{1}{2}\right]\right|\overline{\mbox{rng}(I-A^*A)^\frac{1}{2}}$$
and it is a bounded linear operator between $\mathfrak{D}_A=\overline{\mbox{rng}}(I-A^*A)^\frac{1}{2}$ and $\mathfrak{D}_{A^*}=\overline{\mbox{rng}}(I-AA^*)^\frac{1}{2}$.\\
\vspace{.2cm}

 \textbf{First Situation:}~Suppose $A\in C_{.0}$ is similar to a unitary operator and the spectrum $\sigma(A)$ is a proper subset of $\mathbb{T}$. The existence of nonunitary contractions, specially with a compact spectrum, can be found in \cite{Radj}. Now the characteristic function $C_A(z)$ is boundedly invertible on the open unit disc and is an inner function (see: Theorem $4.5$ in \cite{Nagy}). Also, $C_A(z)$ are unitary operators for every $z$ on the unit circle except $\sigma(A)$. \\
As we have mentioned in section~\ref{Se 3}, we can consider $C_A(z)\in\mathcal{S}$. Also, $C_A(x)$ are unitary operators for all $x\in \mathbb{R}\setminus \mathfrak{S}$, where $\mathfrak{S}$ is the pre-image of $\sigma(A)$ under the conformal map $C$. Now, we can extend $C_A(z)$ to the lower half plane by
$$C_A(z)=\{C_A(\overline{z})^*\}^{-1}\hspace{.5cm}\mbox{for}~z\in\mathbb{C}_-.$$
We denote the extended function as $\mathfrak{C}_A(z)$. Thus we can have an RKHS similar to the one mentioned in lemma \ref{Extension} based on $\mathfrak{C}_A(z)$.
\vspace{.3cm}

\textbf{Second Situation:}~Suppose $A\in C_0$ is a unicellular operator with the scalar multiple equal to the minimal function $m_A(z)$. It is known that the minimal function of this type of operator $A$ is a singular inner function (see: Proposition $7.3$, \cite{Nagy}). Thus $C_A(z)$ is invertible for all $z\in \mathbb{D}$ (by Theorem $5.1$ in \cite{Nagy}). Also, the spectrum $\sigma(A)$ consists of a single point of $\mathbb{T}$, and without loss of any generality, we can assume that $\sigma(A)=\{1\}$. Thus by using the conformal map $C$, we can have an operator valued function $\chi_A(z)$ (say) in $\mathcal{S}^{in}\cap \mathcal{S}^{in}_*$. Moreover, $\chi_A(x)$ is unitary for all $x\in\mathbb{R}$. Now similarly to the first situation, we can extend $\chi_A(z)$ as an operator valued entire function and construct an RKHS $\mathcal{H}(\chi_A)$ of $\mathfrak{X}$-valued entire functions.
\vspace{.3cm}

\textbf{de Branges spaces of entire functions based on a cnu contraction operator:}\\
Let $A$ be a cnu contraction operator as in the second situation, and $E_+$ is a Fredholm operator valued entire function such that $E_+(z)\in B(\mathfrak{D}_{A^*},\mathfrak{D}_A)$ for all $z\in\mathbb{C}$. Also, $E_+(z)$ is invertible at least at one point. Now consider
$$E_-(z)=E_+(z)\chi_A(z)\hspace{.5cm}\mbox{for all}~z\in\mathbb{C}.$$
Thus $E_-$ is a Fredholm operator valued entire function, $E_-(z)\in B(\mathfrak{D}_{A},\mathfrak{D}_A)$ for all $z\in\mathbb{C}$ and $E_-(z)$ is invertible at least at one point. Also
$$E_+^{-1}E_-=\chi_A\in\mathcal{S}^{in}\cap\mathcal{S}^{in}_*.$$
Classes $\mathcal{S}^{in}$ and $\mathcal{S}^{in}_*$ should be understood in the present context. Hence the pair of operator valued functions $(E_-(z),E_+(z))$ for every $z\in\mathbb{C}$ will represent a de Branges operator.

\vspace{.8cm}

\noindent \textbf{Acknowledgements:} 
The authors are thankful to the anonymous referee for pointing out a mistake in Theorem $10.3$ in the early version of this paper and providing several useful suggestions, which eventually led to an improved version.
\vspace{.5cm}

\noindent The authors are also grateful to Professor Harry Dym for 
carefully reading an early version of this paper and suggesting improvements.

\vspace{.5cm}

\noindent The research of the first author is supported by the University Grants Commission (UGC) fellowship (Ref. No. DEC18-424729), Govt. of India.
The research of the second author is supported by the  DST-INSPIRE Faculty research grant (DST/INSPIRE/04/2016/000808) and SERB grant ( SRG/2020/001908 dated 26 October, 2020).

\vspace{.5in}


\end{document}